\documentclass[12pt]{amsart}

\usepackage[latin1]{inputenc}
\usepackage[T1]{fontenc}
\usepackage{amsfonts,amssymb,latexsym}
\usepackage{amsthm}
\usepackage[dvips]{color}

\usepackage{url}

\usepackage{float}
\usepackage{graphicx}

\newtheorem{Th}{Theorem}[section]
\newtheorem{Lem}[Th]{Lemma}

\newtheorem{Prop}[Th]{Proposition}

\newcommand\Z{{\mathbb Z}}
\newcommand\N{{\mathbb N}}
\newcommand\R{{\mathbb R}}

\newcommand\m{{\arrowvert}}
\newcommand\n{{\Arrowvert}}

\begin{document}

\title{Estimation of anisotropic Gaussian fields through Radon transform}

\author{Hermine Bierm\'{e}}
\email{hermine.bierme@math-info.univ-paris5.fr}

\address{MAP5-UMR 8145, Universit\'e Ren\'e Descartes
45, rue des Saints-P\`eres,   75270 PARIS cedex 06  FRANCE,
\url{www.math-info.univ-paris5.fr/~bierme} }
\date{\today}

\author{Fr\'ed\'eric Richard}
\email{frederic.richard@math-info.univ-paris5.fr}
\thanks{This work was supported by ANR grant ``mipomodim'' NT05-1-42030.}

\begin{abstract}
We estimate the anisotropic index of an anisotropic fractional Brownian field. For all directions, we give a convergent
 estimator of the value of the anisotropic index in this direction, based on generalized quadratic variations. We also prove  a central limit theorem. First we present a   result of identification that relies on the asymptotic behavior of the spectral density of a process. Then, we  define Radon transforms of the anisotropic fractional Brownian field and prove that these processes admit a spectral density satisfying the previous assumptions. Finally we use simulated fields to test the proposed estimator in different anisotropic and isotropic cases. Results show that the estimator behaves similarly in  all cases  and  is able to detect anisotropy quite
accurately. 
\end{abstract}

\maketitle

\section*{Introduction}
The one dimensional fractional Brownian motion (fBm) was defined through a stochastic integral by Mandelbrot and Van Ness \cite{MVN} in 1968 for the modeling of
irregular data such as the level of water flows or economic series.
Let us recall that this process is a Gaussian zero mean process with
stationary increments characterized by  its so-called Hurst index
$H\in (0,1)$ and denoted by $B_H=\left\{B_H(t);t\in\R\right\}$. A
generalization of Bochner's Theorem allows to give a spectral
representation of its covariance function, namely
\begin{equation}\label{1DfBm}
\mbox{Cov}\left(B_H(t),B_H(s)\right)=\int_{\R}\left(e^{-it\xi}-1\right)\left(e^{is\xi}-1\right){\m\xi\m^{-2H-1}}d\xi.
\end{equation}
The function ${\m\xi\m^{-2H-1}}$ is called the spectral density of the fBm. Processes with that kind of spectral density are called "$1/f$-noises" in the terminology of signal theory. The Hurst parameter is the index of irregularity of the fBm. It corresponds to  the order of self-similarity of the process and to the critical Hölder exponent of its paths. \\

In this paper we consider  $d$-parameter real-valued Gaussian fields with zero mean and stationary increments,
defined through a spectral representation
\begin{equation}\label{spectral}
\left\{\int_{\R^d}\left(e^{-it\cdot\xi}-1\right)f(\xi)^{1/2}W({d}\xi) ; t\in \R^d\right\},
\end{equation}
where $\cdot$ is the usual scalar product on $\R^d$ and $W$ is a complex Brownian measure with $W(-A)=\overline{W(A)}$ for any Borel set $A\subset \R^d$.
The function  $f$, called the spectral density, is a positive even function of $L^1\left(\R^d,\min\left(1,\m\xi\m^2\right)d\xi\right)$.

A natural extension  of the 1-dimensional fBm is obtained when the spectral density is given by ${\m\xi\m^{-2H-d}}$, where $\m\cdot\m$ is the euclidean norm on $\R^d$. This yields  a zero mean Gaussian field with stationary increments which is isotropic and has therefore the same critical H\"older exponent $H$ in all directions of $\R^d$.

In order  to get an anisotropic field with stationary increments  A. Bonami and A. Estrade define in \cite{ABAE} an anisotropic fractional Brownian field  by considering a spectral density
of the shape  
\begin{equation}\label{density}
f_h(\xi)={\m\xi\m^{-2h(\xi)-d}},
\end{equation} where the power $h(\xi)\in (0,1)$ is an even function which depends on the direction  $\xi/\m\xi\m$ of $\R^d$.
Other generalizations have been proposed for  anisotropic data modeling like the fractional Brownian sheet \cite{Leger} or the multifractional Brownian motion introduced simultaneously in \cite{BJR} and \cite{Peltier}, where the Hurst parameter $H$ is replaced by a  function depending on the point $t\in \R^d$. However such generalizations yield models with non stationary increments.  \\

 Here we consider  an anisotropic fractional Brownian field $X$ defined by \eqref{spectral} with \eqref{density} for some  anisotropic index $h$ and we  focus  on the identification of $h$. As already noticed in \cite{ABAE} this index cannot be recovered by analysing $X$ line by line since its regularity along a line does not depend on the direction. Hence,  to recover anisotropy, authors give another directional analysis method, which is based on field projections. 
Actually, the critical H\"older exponent  of the process $\left\{R_\theta X(t);t\in\R\right\}$  obtained by averaging the field along an hyperplane orthogonal to a fixed direction $\theta$, called Radon transform of $X$, is proved to be equal to  $h(\theta)+\frac{d-1}{2}$ for this direction.
An
estimator of  $h(\theta)$ is then proposed in \cite{ABAEAA}, using quadratic variations to estimate the regularity of the process  $R_\theta X$. However, no speed of convergence
nor asymptotic normality
can be found under their weak conditions  that the spectral density behaves like $f_h$ at high frequencies.\\

Actually, many  estimators for the Hurst parameter of a 1D fBm have been proposed, based for example on time domain methods or spectral methods (see \cite{Coeurjolly} and \cite{BardetE} and references therein).  Quadratic variations can lead to  relevant estimators with asymptotic normality of the Hölder exponent of more general Gaussian processes with stationary increments as proved in \cite{IL} or \cite{Kent} for instance. Moreover in \cite{Lang} the authors give precise bounds of the bias of the variance and show that minimax rates are achieved for this kind of estimators. However,  these previous works  need assumptions on the variogramme of the process $X$ of the following type
\begin{equation}\label{variogramme}
v(t)=\mathbb{E}\left(\left(X(t+t_0)-X(t_0)\right)^2\right)=C|t|^{2H}+r(t) \mbox{ and } r(t)=\underset{t\rightarrow 0}{o}\left(|t|^{2H}\right),
\end{equation}
with further regularity assumptions on the reminder. This leads to a first restriction on the set of values of $H$ since $H$ must be in $(0,1)$. For this reason we adopt here a spectral point of view related to
  the problem of the identification of filtered white noise introduced in \cite{benassiFWN}. In the simplest case of this paper, the authors consider a spectral density $f$ given by
$$f(\xi)=c\m \xi\m^{-2H-1}+R(\xi),$$
with $H\in \R^{+}\smallsetminus \N$, $c>0$ and $R\in{\mathcal C}^2(\R)$, satisfying
$\m R^{(j)}(\xi)\m\le C\m \xi\m^{-2H-1-s-j}$, with $s>0$, for $0\le j\le 2$. It follows that, when $H>1$, one has to consider 
a process no more with stationary increments but with higher order stationary increments.\\ 

In this paper,  under the assumption that the spectral density of a Gaussian process with stationary increments satisfies
for some $c, s>0$
\begin{equation}\label{DASdensite}
f(\xi)= c\m \xi\m^{-2H-1}+\underset{\m\xi\m\rightarrow+\infty}{O}\left(\m\xi\m^{-2H-1-s}\right),
\end{equation} 
we give estimators of $H>0$ using quadratic variations of the process.
On the one hand, for $H\in (0,1)$ this assumption leads
to \eqref{variogramme}, using the fact that $v(t)=4\int_{\R}\sin^2\left(\frac{\xi t}{2}\right)f(\xi)d\xi$.
On the other hand, $H$ can be an integer and we do not need that the reminder is in ${\mathcal C}^2(\R)$ as in \cite{benassiFWN}. \\
Then, with further assumptions on the derivative of $f$ at high frequencies, we get precise error estimates and asymptotic normality. Our main result
is  that the spectral density of the Radon transform of an anisotropic fractional Brownian field satisfies \eqref{DASdensite} with $H=h(\theta)+\frac{d-1}{2}$. Therefore we can get consistent estimators of the anisotropic index $h$  with asymptotic normality using quadratic variations of the Radon transform process.\\

The paper is organized as follows.  The first section is devoted to the estimation of the Hölder exponent of Gaussian processes ($d=1$) with spectral density, using generalized quadratic  variations. We give estimators with asymptotic normality under assumptions that rely on the asymptotic behavior of the spectral density. In the second part we estimate the anisotropic index of an anisotropic fractional Brownian field, using Radon transforms of the field. These transformations lead to  Gaussian processes with spectral densities for which we give an asymptotic expansion. Then we  apply the results of the first part to this process.
In the third part, we test the proposed estimator  using anisotropic and isotropic simulated fields  in two dimensions.

\section{Identification of the exponent for a 1D-process}\label{GQV1D}
We  prove in this section  a first identification result in a general setting. It
is based on the well-known fact that a consistent estimator of the critical  Hölder exponent
of a Gaussian process with stationary increments can be recovered using
generalized quadratic variations (see
\cite{IL} or \cite{Kent} for instance).  Actually, many authors have  considered these estimators
under assumptions based on the variogramme of the process. This is not adapted here for our framework and we 
prove similar results under assumptions based on the
asymptotic behavior of the spectral density. 
Let us recall that, up to a constant, the spectral density of a 1D fractional Brownian motion of Hurst parameter $H \in (0,1)$ is given by the function
${\m \xi\m^{-2H-1}}.$
Remark that no process with stationary increments can admit such spectral density whenever $H\ge 1$ since this function does not belong to $L^1(\R,\min{(1,\m\xi\m^2)}d\xi)$ anymore. However one can obtain such spectral densities by considering processes with higher order stationary increments (see \cite{BJR} and \cite{LESILeana} for instance). Moreover there is no restriction to define a process $X$  with spectral density asymptotically equivalent to that kind of function for any
 $H>0$. Actually, when $H\ge 1$, writing $H=j+s$ with $j\in\N$ and $s\in [0,1)$, $X$ will be $j$ times differentiable in mean square and $X^{(j)}$ will admit $s$ as critical Hölder exponent.\\

Let us consider a zero mean Gaussian process $X=\{X(t);t\in \R \}$, with stationary increments and spectral density
 $f\in L^1(\R,\min{(1,\m \xi\m^{2})}d\xi)$. 
Let us assume that $f$ satisfies \eqref{DASdensite} namely
$$f(\xi)= c\m \xi\m^{-2H-1}+\underset{\m\xi\m\rightarrow+\infty}{O}\left(\m\xi\m^{-2H-1-s}\right),$$
for some $H, c, s>0$. We observe a realization of $X$  at points $\frac{k}{N}$ for $k=0,\ldots,N$:
$$\left\{X(0),X\left(\frac{1}{N}\right),\ldots,X\left(\frac{N}{N}\right)\right\}.$$
Our purpose is  to estimate $H$. The key idea of the former works concerning the estimation of 
H\"older exponent of Gaussian processes with stationary increments as \cite{IL,Kent,Lang,Coeurjolly} for instance, is to consider increments of the process to get a stationary process. For instance, since $X$ has stationary increments the process
$\left\{X\left(\frac{t+1}{N}\right)-X\left(\frac{t}{N}\right); t\in \R\right\}$ is stationary. More generally,
one can consider the filtered process of  $X$
$$Z_{N,\mathbf{a}}(t)=\sum_{k=0}^la_kX\left(\frac{t+k}{N}\right),\,\,\mbox{ for } t\in\R.$$
This defines a stationary process when $\mathbf{a}=(a_0,\ldots,a_l)$ is a discrete filter  of length $l+1$ and of order $K\ge 1$
 ($l, K\in\N$  with $l\ge K$), which means that
$$\sum_{k=0}^la_kk^r=0\,\,\,\,\,\,\mbox{ for } 0\le r\le K-1 \,\,\,\,\,\,\mbox{ and }\,\,\,\,\,\,   \sum_{k=0}^la_kk^{K}\neq 0.$$
For  such a filter $\mathbf{a}$, the observations of $X$ allow us to compute $Z_{N,\mathbf{a}}(p)$ for $p=0,\ldots,N-l$.
Let us remark that if we choose the filter of order $1$ given by $\mathbf{a}=(1,-1)$ then
$Z_{N,\mathbf{a}}(p)=X\left(\frac{p+1}{N}\right)-X\left(\frac{p}{N}\right)$ is the increment of $X$
at point $p/N$ with step $1/N$. More generally, for $K\ge 1$,  the increments of order $K$ of $X$ at point $p/N$ with step $1/N$
are given by $Z_{N,\mathbf{a}}(p)$ for the filter $\mathbf{a}$
of order $K$ and length $K+1$ with $a_k=(-1)^{K-k}\begin{pmatrix} K\\k\end{pmatrix}=(-1)^{K-k}\frac{K!}{k!(K-k)!}$ for $0\le k\le K$.
Let us point out that, since $X$ has stationary increments, 
 $Z_{N,\mathbf{a}}$ is a stationary process for any filter of order $K\ge 1$. 
Following \cite{Kent} we also consider the filtered process of  $X$ with a dilated filter.
More precisely, for an integer $u\ge 1$, the dilation  $\mathbf{a}^u$ of $\mathbf{a}$  is defined 
for $0\le k\le lu$ by
$$a_k^u=\left\{\begin{array}{cc} a_{k'}&\mbox{ if } k=k'u\\
0&\mbox{ otherwise. }\end{array}\right.$$
Since $\underset{k=0}{\overset{lu}{\sum}}k^ra_k^u=u^r\underset{k=0}{\overset{l}{\sum}}k^ra_k$, the filter $\mathbf{a}^u$
has the same order than $\mathbf{a}$. Then, for  a filter $\mathbf{a}$ of length $l$ and of order $K\ge 1$ and $u\ge 1$, due to the stationarity of
the corresponding filtered process $Z_{N,\mathbf{a^u}}$, we can estimate the empirical variance of $Z_{N,\mathbf{a^u}}(0)$ based on the observations of $X$ by considering
\begin{equation} \label{eqVN}
V_{N,\mathbf{a^u}}=\frac{1}{N-lu+1}\sum_{p=0}^{N-lu} \left(Z_{N,\mathbf{a^u}}(p) 
\right) ^2,
\end{equation}
which we call {\it generalized quadratic variations} of $X$. \\

Let us point out that we consider here the same kind of set of locations for the sum as in \cite{Kent,IL} but one could also consider more general one as done in \cite{Lang}.
Moreover, let us remark that one can also consider  $m$-variations of the process  that estimate $\mathbb{E}\left(Z_{N,\mathbf{a^u}}(0)^m\right)$. We will focus here on the quadratic variations ($m=2$). It is motivated by a result of J. F. Coeurjolly \cite{Coeurjolly2} who proves that, in the  fractional Brownian motion  case, the asymptotic variance of the Hurst parameter estimator is the lowest for $m=2$.\\

To build estimators of $H$ the main idea is to choose a filter $\mathbf{a}$ such that  $$ \mathbb{E}\left(V_{N,\mathbf{a^u}}\right)=\mathbb{E}\left(\left(Z_{N,\mathbf{a^u}}(0)\right)^2\right)\underset{N\rightarrow +\infty}{\sim} CN^{-2H}u^{2H}.$$ Then, by considering
estimators given by 
\begin{equation}\label{TN}
T_{N,\mathbf{a^u}}=\frac{V_{N,\mathbf{a^u}}}{\mathbb{E}\left(V_{N,\mathbf{a^u}}\right)},
\end{equation} precise estimates of the asymptotic behavior of
$\mbox{Cov}\left( Z_{N,\mathbf{a^u}}(p),Z_{N,\mathbf{a^v}}(p')\right)$  allow to get 
the almost sure  convergence of $\left(T_{N,\mathbf{a^u}},T_{N,\mathbf{a^v}}\right)$ to $(1,1)$ with asymptotic normality.
Then, an asymptotic estimator of $H$ can be built by considering for instance $\frac{1}{2}\log\left(\frac{V_{N,\mathbf{a^u}}}{V_{N,\mathbf{a^v}}}\right)/\log\left(\frac{u}{v}\right)$ for $u\neq v$.\\

The end of this section is devoted to the rigorous proofs of these statements, under assumptions that rely on the asymptotic
behavior of the spectral density.\\

 Our first result
allows to get an asymptotic development of $\mathbb{E}\left(V_{N,\mathbf{a^u}}\right)$ according to the
asymptotic development of the spectral density at high frequencies.
Let us remark that if we associate to the filter $\mathbf{a}$ the polynomial
$$P_\mathbf{a}(x)=\underset{k=0}{\overset{l}{\sum}}a_kx^k, \mbox{ for } x\in\R,$$ then $\mathbf{a}$ is a filter of order $K$ if and only if 
$P_\mathbf{a}^{(r)}(1)=0$, for $0\le r\le K-1$ and 
$P_\mathbf{a}^{(K)}(1)\neq 0$.

  \begin{Prop}\label{esperance}
  Let  $X=\{X(t);t\in \R\}$ be a zero mean Gaussian process,
with stationary increments and spectral density
 $f$.
We assume that  that there exists $H>~0$, $c>0$ and $s> 0$ such that $f$ satisfies \eqref{DASdensite}.
Then for any filter $\mathbf{a}$ of order $K\ge 1$ and $u\ge 1$ we have the following asymptotics:
$$N^{2K}\mathbb{E}\left(V_{N,\mathbf{a^u}}\right)=\left\{\begin{array}{lll}
u^{2K}\left(\frac{P_{\mathbf{a}}^{(K)}(1)}{K!}\right)^2\int_{\R}\xi^{2K}f(\xi)d\xi&+\underset{N\rightarrow +\infty}{o}
\left(1\right) &\mbox{ if } K<H\\
u^{2K}\left(\frac{P_{\mathbf{a}}^{(K)}(1)}{K!}\right)^22c\log N &+\underset{N\rightarrow +\infty}{o}
\left(\log N\right) &\mbox{ if } K=H.\end{array}\right.$$
Moreover, if  K>H then
$E_{\mathbf{a}}^u(H)=u^{2H}\int_{\R}\left|P_{\mathbf{a}}(e^{-i\xi})\right|^2|\xi|^{-2H-1}d\xi<+\infty$ and
$$N^{2H}\mathbb{E}\left(V_{N,\mathbf{a^u}}\right)=cu^{2H}E_{\mathbf{a}}^1(H)
+\left\{
\begin{array}{ll}
\underset{N\rightarrow +\infty}{O}\left(N^{-2(K-H)}\right) &\mbox{ if } K-H<s/2\\
\underset{N\rightarrow +\infty}{O}\left(N^{-s}\log N\right) &\mbox{ if } K-H=s/2\\
\underset{N\rightarrow +\infty}{O}\left(N^{-s}\right) &\mbox{ if } K-H>s/2
\end{array}\right.$$
  \end{Prop}
\begin{proof}
Since by assumption  $X$ has stationary increments and spectral density $f$ it follows  from (\ref{spectral}) for $d=1$ that, when $K\ge 1$,
$$
Z_{N,\mathbf{a^u}}(t)\stackrel{L^2(\Omega)}{=}\int_{\R}e^{-i\frac{t\xi}{N}}P_{\mathbf{a}}\left(e^{-i\frac{u\xi}{N}}\right)
f(\xi)^{1/2}W(d\xi).
$$
Therefore, for all $p\in \Z$,
$$\mathbb{E}\left(\left(Z_{N,\mathbf{a^u}}(p)\right)^2\right)=\int_{\R}\left|P_{\mathbf{a}}\left(e^{-i\frac{u\xi}{N}}\right)\right|^2 f(\xi)d\xi=\mathbb{E}\left(V_{N,\mathbf{a^u}}\right).$$
Let us assume that $K\le H$, since $\mathbf{a}$ is of order $K$, from Taylor formula, for $\epsilon>0$ there exists
$\delta>0$ such that $\left|\left|P_{\mathbf{a}}\left(e^{-i{\xi}}\right)\right|^2-\left(\frac{P_{\mathbf{a}}^{(K)}(1)}{K!}\right)^2\xi^{2K}\right|\le
\epsilon\xi^{2K}.$
Hence,
$$\left|
\int_{|\xi|\le \frac{\delta N}{u}}\left|P_{\mathbf{a}}\left(e^{-i\frac{u\xi}{N}}\right)\right|^2 f(\xi)d\xi
-N^{-2K}u^{2K}\left(\frac{P_{\mathbf{a}}^{(K)}(1)}{K!}\right)^2\int_{|\xi|\le \frac{\delta N}{u}}\xi^{2K}f(\xi)d\xi\right|$$
$$\le \epsilon N^{-2K}\int_{|\xi|\le \frac{\delta N}{u}}\xi^{2K}f(\xi)d\xi.$$
According to \eqref{DASdensite},
$$\int_{|\xi|\le \frac{\delta N}{u}}\xi^{2K}f(\xi)d\xi=\left\{\begin{array}{lll}\int_{\R}\xi^{2K}f(\xi)d\xi&+
\underset{N\rightarrow +\infty}{o}(1)&\mbox{if } K<H\\
2c\log N &+
\underset{N\rightarrow +\infty}{O}(1)&\mbox{if } K=H
\end{array}\right.$$
and $\int_{|\xi|>\frac{\delta N}{u}}\left|P_{\mathbf{a}}\left(e^{-i\frac{u\xi}{N}}\right)\right|^2 f(\xi)d\xi=\underset{N\rightarrow +\infty}{O}(N^{-2H})$, which yields the result.\\

Now, let us assume that $K>H$ and write $f( \xi)=c|\xi|^{-2H-1}+R(\xi)$ with $R\in~L^1\left(\R,\min(1,|\xi|^{2K})d\xi\right)$
satisfying $R(\xi)=\underset{\m \xi\m\rightarrow +\infty}{O}\left(\m \xi\m^{-2H-1-s}\right)$. A change of variables leads to
$$\int_{\R}\left|P_{\mathbf{a}}\left(e^{-i\frac{u\xi}{N}}\right)\right|^2 |\xi|^{-2H-1}d\xi=N^{-2H}u^{2H}E_{\mathbf{a}}^1(H),$$
with $E_{\mathbf{a}}^1(H)=\int_{\R}\left|P_{\mathbf{a}}\left(e^{-i\xi}\right)\right|^2 |\xi|^{-2H-1}d\xi<\infty$ since $H<K$. Then, applying the previous results to $R$ with $H$ replaced by $H+s/2$ we obtain the result for the reminder term, which concludes the proof.
\end{proof}

Therefore, to recover the parameter $H$, one has to consider a filter $\mathbf{a}$ of order $K>H$. In this case, for any $u\ge 1$, by stationarity of $Z_{N,\mathbf{a^u}}$ we obtain that 
$\mathbb{E}\left(V_{N,\mathbf{a^u}}\right)\underset{N\rightarrow +\infty}{\sim} N^{-2H}u^{2H}cE_{\mathbf{a}}^1(H)$. In order to prove the almost sure  convergence of $\left(T_{N,\mathbf{a^u}},T_{N,\mathbf{a^v}}\right)$ for $u,v\ge 1$, with $T_{N,\mathbf{a^u}}$ given by \eqref{TN}, one has to estimate
$\mbox{Cov}\left(V_{N,\mathbf{a^u}},V_{N,\mathbf{a^v}}\right)$.
Actually,
$$\mbox{Cov}\left(V_{N,\mathbf{a^u}},V_{N,\mathbf{a^v}}\right)=\frac{1}{N-lu+1}\frac{1}{N-lv+1}\sum_{p=0}^{N-lu}\sum_{p=0}^{N-lv} \mbox{Cov}\left( Z_{N,\mathbf{a^u}}(p)^2,Z_{N,\mathbf{a^v}}(p')^2\right),$$
with  $\mbox{Cov}\left( Z_{N,\mathbf{a^u}}(p)^2,Z_{N,\mathbf{a^v}}(p')^2\right)=2\mbox{Cov}\left(Z_{N,\mathbf{a^u}}(p),Z_{N,\mathbf{a^v}}(p')\right)^2$
since $(Z_{N,\mathbf{a^u}}(p),Z_{N,\mathbf{a^v}}(p'))$ is a Gaussian vector.
Moreover, 
$$\mbox{Cov}\left(Z_{N,\mathbf{a^u}}(p),Z_{N,\mathbf{a^v}}(p')\right)=\int_{\R}e^{-i\frac{(p-p')\xi}{N}}P_{\mathbf{a}}\left(e^{-i\frac{u\xi}{N}}\right)\overline{P_{\mathbf{a}}\left(e^{-i\frac{v\xi}{N}}\right)}f(\xi)d\xi.$$
Let us denote 
\begin{equation}\label{gammah}
\Gamma_{N,\mathbf{a}}^{u,v}(p)=\int_{\R}e^{-i\frac{p\xi}{N}}h_{\mathbf{a}}^{u,v}\left(\frac{\xi}{N}\right)f(\xi)d\xi \mbox{ where }
h_{\mathbf{a}}^{u,v}(\xi)=P_{\mathbf{a}}\left(e^{-iu\xi}\right)\overline{P_{\mathbf{a}}\left(e^{-iv\xi}\right)}
\end{equation}
such that
$$\mbox{Cov}\left(V_{N,\mathbf{a^u}},V_{N,\mathbf{a^v}}\right)=\frac{2}{N-lu+1}\sum_{p=-N+lv}^{N-lu}\Gamma_{N,\mathbf{a}}^{u,v}(p)^2.$$
It is obvious that under \eqref{DASdensite}, for $K>H$ we have
$$\Gamma_{N,\mathbf{a}}^{u,v}(p)\underset{N\rightarrow +\infty}{\sim}cN^{-2H}\int_{\R}e^{-ip\xi}h_{\mathbf{a}}^{u,v}\left(\xi\right)|\xi|^{-2H-1}d\xi.$$
However, we have to consider  $\underset{p=-N+lv}{\overset{N-lu}{\sum}} \Gamma_{N,\mathbf{a}}^{u,v}(p)^2$ and further assumptions have to be done to see when
$\left(\int_{\R}e^{-ip\xi}h_{\mathbf{a}}^{u,v}\left(\xi\right)|\xi|^{-2H-1}d\xi\right)_{p\in\Z}$ is in $\ell^2(\Z)$.

\begin{Prop}\label{variance} Under the assumptions of Proposition \ref{esperance}, if we assume moreover that
$f$ is differentiable on $\R\smallsetminus(-r,r)$, for r large enough
and
\begin{equation}\label{DASD}
f'(\xi)=-(2H+1)\frac{c}{\m \xi\m^{2H+2}}+o_{\m \xi\m\rightarrow +\infty}\left(\frac{1}{\m \xi\m^{2H+2}}\right),
\end{equation}
then, for $K>H$ and any $\delta<\min(2(K-H),1)$  with $\delta>\max(1-2H,0)$, one can find $C>0$ such that, for all $|p|\le N$,
$$\left|\Gamma_{N,\mathbf{a}}^{u,v}(p)\right|\le CN^{-2H}\left(1+|p|\right)^{-\delta}.$$
Moreover, for any $K>H+1/4$ we have $C_{\mathbf{a}}^{u,v}(H)=2\underset{p\in\Z}{\sum}\left(\int_{\R}e^{-ip\xi}h_{\mathbf{a}}^{u,v}\left(\xi\right)|\xi|^{-2H-1}d\xi\right)^2~<~+\infty$
and
$$\mbox{Cov}\left(V_{N,\mathbf{a^u}},V_{N,\mathbf{a^v}}\right)\underset{N\rightarrow +\infty}{\sim}
c^2C_{\mathbf{a}}^{u,v}(H)N^{-4H-1}.$$
\end{Prop}
\begin{proof}
Let $K>H$.
Let us write $f(\xi)=c|\xi|^{-2H-1}+R(\xi)$ with $R\in L^1\left(\R,\min(1,|\xi|^{2K})d\xi\right)$. 
By a change of variables, we can write 
$$\Gamma_{N,\mathbf{a}}^{u,v}(p)=cN^{-2H}\int_{\R}e^{-ip\xi}h_{\mathbf{a}}^{u,v}\left(\xi\right)|\xi|^{-2H-1}d\xi
+\int_{\R}e^{-i\frac{p\xi}{N}}h_{\mathbf{a}}^{u,v}\left(\frac{\xi}{N}\right)R(\xi)d\xi.$$
Let us focus on the first term and remark that for $u=v$, since $h_{\mathbf{a}}^{u,u}(\xi)=\left|P_{\mathbf{a}}\left(e^{-iu\xi}\right)\right|^2$ is real 
$$\int_{\R}e^{-ip\xi}h_{\mathbf{a}}^{u,u}\left(\xi\right)|\xi|^{-2H-1}d\xi=\int_{\R}\cos(p\xi)h_{\mathbf{a}}^{u,u}\left(\xi\right)|\xi|^{-2H-1}d\xi.$$
Moreover, since $\mathbf{a}$ is a filter of order $K$, by Taylor formula one can find $C>0$ such that
\begin{equation}\label{ha}
\left|h_{\mathbf{a}}^{u,u}\left(\xi\right)\right|\le C\min\left(1,\xi^{2K}\right) \mbox{ and }
\left|\frac{\rm d}{{\rm d}\xi}h_{\mathbf{a}}^{u,u}\left(\xi\right)\right|\le C\min\left(1,|\xi|^{2K-1}\right).
\end{equation}
Therefore, when $p\neq 0$,
we can integrate by parts
$$\int_{\R}e^{-ip\xi}h_{\mathbf{a}}^{u,u}\left(\xi\right)|\xi|^{-2H-1}d\xi=-\int_{\R}\frac{\sin(p\xi)}{p}\frac{\rm d}{{\rm d}\xi}\left(h_{\mathbf{a}}^{u,u}\left(\xi\right)|\xi|^{-2H-1}\right)d\xi.$$
Then, for any $\delta<\min(2(K-H),1)$ with $\delta>\max(1-2H,0)$,
$$\left|\int_{\R}e^{-ip\xi}h_{\mathbf{a}}^{u,u}\left(\xi\right)|\xi|^{-2H-1}d\xi\right|\le |p|^{-\delta}\left|
\int_{\R}|\xi|^{1-\delta}\left|\frac{\rm d}{{\rm d}\xi}\left(h_{\mathbf{a}}^{u,u}\left(\xi\right)|\xi|^{-2H-1}\right)\right|d\xi\right|,$$
with $\int_{\R}|\xi|^{1-\delta}\left|\frac{\rm d}{{\rm d}\xi}\left(h_{\mathbf{a}}^{u,u}\left(\xi\right)|\xi|^{-2H-1}\right)\right|d\xi<+\infty$ according to \eqref{ha}.
Writing $h_{\mathbf{a}}^{u,v}\left(\xi\right)$ as
$$\frac{1}{2}\left(\left|P_{\mathbf{a}}\left(e^{-iu\xi}\right)+P_{\mathbf{a}}\left(e^{-iv\xi}\right)\right|^2+i\left|P_{\mathbf{a}}\left(e^{-iu\xi}\right)+iP_{\mathbf{a}}\left(e^{-iv\xi}\right)\right|^2-(1+i)\left(h_{\mathbf{a}}^{u,u}+h_{\mathbf{a}}^{v,v}\right)\right),$$
we also get that for any $\delta<\min(2(K-H),1)$ with $\delta>\max(1-2H,0)$, for all $p\in\N$,
$$\left|\int_{\R}e^{-ip\xi}h_{\mathbf{a}}^{u,v}\left(\xi\right)|\xi|^{-2H-1}d\xi\right|\le C(1+|p|)^{-\delta}.$$
Therefore $\left(\int_{\R}e^{-ip\xi}h_{\mathbf{a}}^{u,v}\left(\xi\right)|\xi|^{-2H-1}d\xi\right)_{p\in\Z}$ is in $\ell^2(\Z)$ as soon as one can choose $\delta>1/2$, which is possible when $K>H+\frac{1}{4}$.
It remains to study the second term.
Let $\epsilon>0$ and choose $r>0$ such that
$R$ is differentiable on $\R\smallsetminus(-r,r)$ with $|R^{(j)}(\xi)|\le \epsilon |\xi|^{-2H-1-j}$ for $j\in\{0,1\}$ and 
$|\xi|\ge r$. We write
$$\int_{\R}e^{-i\frac{p\xi}{N}}h_{\mathbf{a}}^{u,v}\left(\frac{\xi}{N}\right)R(\xi)d\xi=\int_{|\xi|<r}+\int_{|\xi|\ge r}.$$
For the first integral, since
$R\in L^1(\R,\min{(1,\m\xi\m^{2K})}d\xi)$, we remark that 
$$\left|\int_{|\xi|<r}e^{-i\frac{p\xi}{N}}h_{\mathbf{a}}^{u,v}\left(\frac{\xi}{N}\right)R(\xi)d\xi \right|\le C(r)N^{-2K},$$
where $C(r)$ is a constant that only depends on $r$, which may change line by line.
For the second integral, for all $p \in \Z$, similar computations as previously leads to
$$\left|\int_{|\xi|\ge r}e^{-i\frac{p\xi}{N}}h_{\mathbf{a}}^{u,v}\left(\frac{\xi}{N}\right)R(\xi)d\xi \right|\le C(r)N^{-2K}
+C\epsilon (1+|p|)^{-\delta}N^{-2H},$$
 for any $\delta<\min(2(K-H),1)$ with $\delta>\max(1-2H,0)$. Therefore, for any $|p|\le N$,
 $$\left|\int_{\R}e^{-i\frac{p\xi}{N}}h_{\mathbf{a}}^{u,v}\left(\frac{\xi}{N}\right)R(\xi)d\xi \right|
 \le \left(C(r)N^{-2(K-H)-\delta}+C\epsilon\right) (1+|p|)^{-\delta}N^{-2H}.$$
 It follows that  one can find $C>0$ such that, for all $|p|\le N$,
$$\left|\Gamma_{N,\mathbf{a}}^{u,v}(p)\right|\le CN^{-2H}\left(1+|p|\right)^{-\delta} \mbox{ with } N^{2H}\Gamma_{N,\mathbf{a}}^{u,v}(p)\underset{N\rightarrow +\infty}{\longrightarrow}c\int_{\R}e^{-ip\xi}h_{\mathbf{a}}^{u,v}\left(\xi\right)|\xi|^{-2H-1}d\xi.$$
Finally, when $K>H+1/4$ one can choose $\delta \in (1/2,1)$ such that
$C_{\mathbf{a}}^{u,v}(H)=2\sum_{p\in\Z}\left(\int_{\R}e^{-ip\xi}h_{\mathbf{a}}^{u,v}\left(\xi\right)|\xi|^{-2H-1}d\xi\right)^2$ is finite. By the dominated convergence theorem,
$$
N^{4H+1}\mbox{Cov}\left(V_{N,\mathbf{a^u}},V_{N,\mathbf{a^v}}\right)=\frac{2N}{N-lu+1}\sum_{p=-N+lv}^{N-lu} \left(N^{2H}\Gamma_N^{u,v}(p)\right)^2\underset{N\rightarrow +\infty}{\longrightarrow} c^2C_{\mathbf{a}}^{u,v}(H).$$

 \end{proof}

We can now state our first identification result.
\begin{Prop} \label{propVNDAS}
Let  $X=\{X(t);t\in \R\}$ be a zero mean Gaussian process,
with stationary increments and spectral density
 $f$, which satisfies the assumptions of Proposition \ref{esperance} and Proposition \ref{variance}.
Let $\mathbf{a}$ be a filter of order $K>H$ and $u,v\ge 1$ two integers with $u\neq v$. Then, almost surely,
$$\widehat{H_{N,\mathbf{a}}}(u,v)=\frac{1}{2\log(u/v)}\log\left(\frac{V_{N,\mathbf{a^u}}}{V_{N,\mathbf{a^v}}}\right)\underset{N\rightarrow+\infty}\longrightarrow H.$$
Moreover, for $K>H+1/4$, let us denote
$$\gamma_\mathbf{a}^{u,v}(H)=\frac{1}{4\log^2(u/v)}\left(\frac{C_\mathbf{a}^{u,u}(H)}{E_\mathbf{a}^u(H)^2}+\frac{C_\mathbf{a}^{v,v}(H)}{E_\mathbf{a}^v(H)^2}
-2\frac{C_\mathbf{a}^{u,v}(H)}{E_\mathbf{a}^u(H)E_\mathbf{a}^v(H)}\right).$$
Then,  when  $s>\frac{1}{2}$,  
$\sqrt{N}\left(\widehat{H_{N,\mathbf{a}}}(u,v)-H\right)\overset{d}{\underset{N\rightarrow +\infty}{\longrightarrow}}{\mathcal N}\left(0,\gamma_\mathbf{a}^{u,v}(H)\right),$ with
$$N\mathbb{E}\left(\left(\widehat{H_{N,\mathbf{a}}}(u,v)-H\right)^2\right)\underset{N\rightarrow +\infty}{\longrightarrow} \gamma_\mathbf{a}^{u,v}(H)$$
 and, when $s\le \frac{1}{2}$, 
$$ \mathbb{E}\left(\left(\widehat{H_{N,\mathbf{a}}}(u,v)-H\right)^2\right)=\underset{N\rightarrow +\infty}{O}\left(N^{-2s}\right).$$
\end{Prop}
\begin{proof} 
Following classical computations on Gaussian quadratic forms as in \cite{IL} for instance, from Proposition \ref{esperance}
and \ref{variance} we obtain that for all $K>H$
$$\left(T_{N,\mathbf{a^u}},T_{N,\mathbf{a^v}}\right)=\left(\frac{V_{N,\mathbf{a^u}}}{\mathbb{E}\left(V_{N,\mathbf{a^u}}\right)},
\frac{V_{N,\mathbf{a^v}}}{\mathbb{E}\left(V_{N,\mathbf{a^v}}\right)}\right)\underset{N\rightarrow +\infty}{\longrightarrow}
(1,1) \mbox{ a.s. }$$
with for all $K>H+1/4$,
$\sqrt{N}\left(T_{N,\mathbf{a^u}},T_{N,\mathbf{a^v}}\right)\overset{d}{\underset{N\rightarrow +\infty}{\longrightarrow}}
{\mathcal N}\left(0,\Sigma_\mathbf{a}^{u,v}(H)\right)$
\mbox{ and }$$ N\mbox{Cov}\left(\left(T_{N,\mathbf{a^u}},T_{N,\mathbf{a^v}}\right)^t\left(T_{N,\mathbf{a^u}},T_{N,\mathbf{a^v}}\right)\right)
\underset{N\rightarrow +\infty}{\longrightarrow}\Sigma_\mathbf{a}^{u,v}(H),$$
where $$\Sigma_\mathbf{a}^{u,v}(H)=\left(\begin{array}{cc} {C_\mathbf{a}^{u,u}(H)}/{E_\mathbf{a}^u(H)^2}& {C_\mathbf{a}^{u,v}(H)}/{E_\mathbf{a}^u(H)E_\mathbf{a}^v(H)}\\
{C_\mathbf{a}^{u,v}(H)}/{E_\mathbf{a}^u(H)E_\mathbf{a}^v(H)}&{C_\mathbf{a}^{v,v}(H)}/{E_\mathbf{a}^v(H)^2}\end{array}\right).$$
Then, using Taylor Formula for the function $g(x,y)=\log\left(\frac{x}{y}\right)$ (see Theorem 3.3.11 in \cite{dacunha} for instance)
we get that for  $K>H$ almost surely
$\log\left(\frac{T_{N,\mathbf{a^u}}}{T_{N,\mathbf{a^v}}}\right)\underset{N\rightarrow +\infty}{\longrightarrow} 0 $
with, for all $K>H+1/4$,
$\sqrt{N}\log\left(\frac{T_{N,\mathbf{a^u}}}{T_{N,\mathbf{a^v}}}\right)
\overset{d}{\underset{N\rightarrow +\infty}{\longrightarrow}}
{\mathcal N}\left(0,\Gamma_\mathbf{a}^{u,v}(H)\right)$
 and $$N\mathbb{E}\left(\log\left(\frac{T_{N,\mathbf{a^u}}}{T_{N,\mathbf{a^v}}}\right)^2\right)
\underset{N\rightarrow +\infty}{\longrightarrow}\Gamma_\mathbf{a}^{u,v}(H),$$
where $\Gamma_\mathbf{a}^{u,v}(H)=4\log(u/v)^2\gamma_\mathbf{a}^{u,v}(H).$
Then let us write for $u\neq v$
$$\widehat{H_{N,\mathbf{a}}}(u,v)=\frac{1}{2\log(u/v)}\left(\log\left(\frac{\mathbb{E}\left(V_{N,\mathbf{a^u}}\right)}{\mathbb{E}\left(V_{N,\mathbf{a^v}}\right)}\right)+\log\left(\frac{T_{N,\mathbf{a^u}}}{T_{N,\mathbf{a^v}}}\right)\right).$$
According to Proposition \ref{esperance}
it is straightforward to see that for $K>H$
$$\widehat{H_{N,\mathbf{a}}}(u,v)\underset{N\rightarrow +\infty}{\longrightarrow} H \mbox{ a.s.}$$
Moreover, for $K>H+1/4$ and $s>1/2$, Proposition \ref{esperance} leads to 
$$\frac{\mathbb{E}\left(V_{N,\mathbf{a^u}}\right)}{\mathbb{E}\left(V_{N,\mathbf{a^v}}\right)}=\left(\frac{u}{v}\right)^{2H}\left(1+\underset{N\rightarrow +\infty}{o}(1/\sqrt{N})\right).$$ So in this case  
$\sqrt{N}\left(\widehat{H_{N,\mathbf{a}}}(u,v)-H\right)\overset{d}{\underset{N\rightarrow +\infty}{\longrightarrow}}
{\mathcal N}\left(0,\gamma_\mathbf{a}^{u,v}(H)\right)$\mbox{ and } $$N\mathbb{E}\left(
\left(\widehat{H_{N,\mathbf{a}}}(u,v)-H\right)^2\right)\underset{N\rightarrow +\infty}{\longrightarrow}\gamma_\mathbf{a}^{u,v}(H).$$
Let us point out that if $s\le 1/2$, from Proposition \ref{esperance} using the fact that $K>H+1/4>H+s/2$ we get 
$$\frac{\mathbb{E}\left(V_{N,\mathbf{a^u}}\right)}{\mathbb{E}\left(V_{N,\mathbf{a^v}}\right)}=\left(\frac{u}{v}\right)^{2H}\left(1+\underset{N\rightarrow +\infty}{O}({N^{-s}})\right).$$
Therefore, 
$$\mathbb{E}\left(
\left(\widehat{H_{N,\mathbf{a}}}(u,v)-H\right)^2\right)=\underset{N\rightarrow +\infty}{O}\left({N^{-2s}}\right).$$
\end{proof}

This estimator is used in the next section in order to estimate the anisotropic index of an anisotropic fractional Brownian
field.
\section{Identification of the anisotropic index of anisotropic fractional Brownian fields}
In this section we consider an anisotropic fractional Brownian field $$X=\left\{X(t) \,;\,\, t\in\R^d\right\},$$ as introduced in \cite{ABAE}, which is a zero mean Gaussian random field, with stationary increments and spectral representation \eqref{spectral}. The
 spectral density is given by \eqref{density}, namely $f_h(\xi)={\m\xi\m^{-2h(\xi)-d}}$,  
where $h$ is an even homogeneous function of degree $0$ with values in
$(0,1)$, called anisotropic index. To determine anisotropy of such a
field one could try to estimate its directional regularity by
extracting lines of the field along various directions. However, for
anisotropic fractional Brownian fields, this method fails. Actually,
when $\theta$ is a fixed direction of the sphere $S^{d-1}$, one can
prove that  the process $\left\{X(t\theta); t\in\R\right\}$ is
still a zero mean Gaussian process with
spectral density given by $$p\in\mathbb{R}\mapsto \int_{\langle \theta \rangle^{\perp}}f_h\left(p\theta+\gamma\right)d\gamma
=\int_{\langle \theta \rangle^{\perp}}p^{-2h(\theta+\gamma)-1}\left(1+|\gamma|^2\right)^{-h(\theta+\gamma)-d/2}d\gamma,$$
where $\langle\theta\rangle^{\perp}$ stands for the hyperplane orthogonal to
$\theta$.
Therefore, according to Proposition 3.6 and Proposition 3.3 of \cite{ABAE}
this process admits a critical Hölder exponent equals to the essential
infimum $h_0$ of the function $h$.
 Then, the study of the generalized quadratic variations of such a process can at most allow us to recover $h_0$. To deal with this obstruction and in order to study processes rather than fields, the authors of \cite{ABAE} have introduced the Radon transform of such fields.\\
When a function $f$ is integrable over $\R^d$, one can define its
Radon transform on $\R$ (see \cite{ramm} for instance), in the
direction $\theta$,  by
$${\mathcal R}_{\theta}f(t)=\int_{\langle\theta\rangle^{\perp}}f(s+t\theta)ds, \mbox{ for all } t\in \R.$$
 For a function $f$ which does not decay sufficiently at
infinity, one can  integrate it against a window. Let
$\rho$ be a smooth function defined on $\langle\theta\rangle^{\perp}$,  that
compensates the behavior at infinity of $f$. Then one can define the
windowed Radon transform of $f$ on $\R$, in the direction $\theta$,
by
$${\mathcal R}_{\theta,\rho}f(t)=\int_{\langle\theta\rangle^{\perp}}f(s+t\theta)\rho(s)ds, \mbox{ for all } t\in \R.$$
We should use strong assumptions on the anisotropic fractional Brownian field $X$ in order to define its windowed Radon transform as an integral. However, 
according to Proposition 4.1 of \cite{ABAE}, one can define the
Radon transform of $X$, with a convenient  window $\rho$, in the direction
$\theta$, by a discretization of the integral. For notational sake of simplicity
we deal with the direction $\theta$  to be $\theta_0=(0,\ldots,0,1)$ and identify the space $\R^{d-1}\times\{0\}$ to $\R^{d-1}$.
Let us choose $\rho$ a  function of the Schwartz class ${\mathcal S}\left(\R^{d-1}\right)$, with real values, ie $\rho$ is a smooth  function rapidly decreasing
\begin{equation}\label{assonthewindow}
\forall N\in\N,\,\,\forall x \in \R^{d-1},\,\,\,\m \rho(x) \m \le C_N(1+\m x \m)^{-N}.
\end{equation}
Then, the process
$$ 2^{-n(d-1)}\sum_{s\in 2^{-n}\Z^{d-1}}
 X\left(s,t\right) \rho(s),\,\,\forall t\in \R,$$
 admits a limit in $L^2(\Omega)$ for  the finite
dimensional distributions, when $n$ tends to infinity.
This limit is called the Radon transform of $X$ with the window $\rho$  and is
denoted by ${ R}_{\rho}X=\{{R}_{\rho}X(t); t\in\R\}$.\\

Let us remark that one can define  the
Radon transform of $X$, with the window $\rho$, under less restrictive assumptions on $\rho$, as soon as the previous limit exists.
  The existence of the limit process is proved  in both
\cite{ABAE} and \cite{HB} and relies on the slow increase of the
covariance function of $X$ due to its stationary increments and on
its mean square continuity. \\
Let us also point out that one can define the Radon transform of $X$ with the window $\rho$ for any direction $\theta$.
Actually, it is sufficient to chose $\kappa_{\theta}$ a rotation of $\R^d$ that maps $\theta_0=(0,\ldots,0,1)$ onto $\theta$. Since $X\circ\kappa_{\theta}$ is still an anisotropic fractional Brownian field,
with anisotropic index given by $h\circ\kappa_{\theta}$, which satisfies the same assumption as $h$, we just have to consider the Radon transform of this field 
with the window $\rho$.\\

 By linearity of such a transformation, the Radon transform of $X$
with the window $\rho$ is still a zero mean Gaussian process with stationary
increments and it admits a spectral density given by the Radon
transform of the spectral density $f_h$ of $X$, given by (\ref{density}), against the window
$\left|\widehat{\rho}\right|^2$,
\begin{equation}\label{TRDS}
{\mathcal R}_{\left|\widehat{\rho}\right|^2}f_h(p)=\int_{\R^{d-1}}f_h(\gamma,p)\left|\widehat{\rho}(\gamma)\right|^2d\gamma,\mbox{ for all } p\in \R,
\end{equation}
where $\widehat{\rho}$ is the $(d-1)$-dimensional Fourier
transform of the window $\rho$. To estimate the Hölder regularity of
this process, we will use its generalized quadratic variations as
introduced in the section \ref{GQV1D}. Therefore, we have  to study the
asymptotic behavior of ${\mathcal
R}_{\left|\widehat{\rho}\right|^2}f_h$ in order to apply Proposition
\ref{propVNDAS}. We  prove and use the following general result
on the windowed Radon transform.
\begin{Prop}\label{DASRT}
Let $h$ and $c$ be  given functions on $\R^d$. Let   $\alpha>0$. We assume that
 $h$ and $c$ are even homogeneous functions of degree  $0$, Lipschitz of order $\alpha$ on the sphere, with $h$ positive.\\
Let $\delta_0>0$. Let  $f$ be a function defined on $\R^d$  such that,  for all $\delta\in (0,\delta_0)$,
$$f(\xi)= \frac{c(\xi)}{ {\m \xi \m}^{h(\xi)}}+ {o}\left(\frac{1}{\m \xi
\m^{h(\xi)+\delta}}\right) \mbox{ when }
{\m \xi \m \rightarrow +\infty}.$$
Choose  $\rho\in
{\mathcal S}\left(\R^{d-1}\right)$ such that $\int_{\R^{d-1}}
\rho(\gamma)d\gamma=1$.
Then, the Radon transform of $f$ with the window $\rho$ satisfies, for all $\delta \in (0,\delta_1)$,
$${\mathcal R}_{\rho}f(p) =\frac{c(\theta_0)}{ {\m  p \m}^{h(\theta_0)}}+
o\left(\frac{1}{\m p \m^{h(\theta_0)+\delta}}\right)\mbox{ when } p\in \R \mbox{ and }
{\m p \m \rightarrow +\infty},$$
with $\delta_1=\min\left(\delta_0, \alpha\right).$
\end{Prop}
\begin{proof}
Let   $\rho$ be a function of ${\mathcal S}\left(\R^{d-1}\right)$ with $\int_{\R^{d-1}}
\rho(\gamma)d\gamma=1$. For
$p \in \R$, with $\m p\m$ large enough, one can define the integral  $${\mathcal R}_{\rho}f( p)=\int_{\R^{d-1}}f(\gamma,p)
 \rho(\gamma)d\gamma.$$
We want to estimate its asymptotics when $\m p \m \rightarrow
+\infty$. First, let us assume that there exists $A>1$ such that,
for $\xi\in\R^d$ and $\m\xi\m>A$,
$$f(\xi)=\frac{c(\xi)}{ {\m \xi \m}^{h(\xi)}},$$
with $h$ and $c$ satisfying assumptions of Proposition \ref{DASRT}.
In this case, we will prove that  for all
$0<\delta<\alpha$,
\begin{equation}\label{TRAFBc}
{\mathcal R}_{\rho}f(p)=f(p\theta_0)+o(\m p \m^{-h(\theta_0)-\delta})\mbox{ when }p\in \R \mbox{ and }
{\m p \m \rightarrow +\infty}.
\end{equation}
For $\m p\m>A$, since $\int_{\R^{d-1}}\rho(\gamma)d\gamma=1$, let us write
$${\mathcal R}_{\rho}f(p)=f(p\theta_0)+\int_{\R^{d-1}}\left(f(\gamma,p)
-f(p\theta_0)\right)\rho(\gamma)d\gamma.$$ Then, it is enough to
give an upper bound for
$$\int_{\R^{d-1}}\left(f(\gamma,p)
-f(p\theta_0)\right)\rho(\gamma)d\gamma.$$
Let us denote $h\left(S^{d-1}\right)=[h_0,h_1]$ with $h_0>0$ by assumption on $h$.
Since $\rho$ is rapidly decreasing, for all $s>0$ and $N\in\N$,
$$\int_{\m \gamma \m>\m p \m^{s}}\left(f(\gamma,p)-f(p\theta_0)\right)
 \rho(\gamma)d\gamma={O}_{\m p \m \rightarrow +\infty}(\m p \m^{-h_0-Ns}),$$
which is negligible compared to $\m p \m^{-h(\theta_0)-\delta}$ as soon as $N>\frac{\delta+h_1-h_0}{s}$.\\
\\
Thus, it is sufficient to consider
$$\Delta_s(p)=\int_{\m \gamma \m \le \m p \m^{s}}\left(f(\gamma,p)
-f(p\theta_0)\right)\rho(\gamma)d\gamma.$$
But,
\begin{eqnarray*}
\left|\Delta_s(p)\right| &\le & \int_{\m \gamma \m \le \m p \m^{s}}\m c(\gamma,p) \m\left|
\frac{1}{ {\left(\m\gamma\m^2 +p^2\right)}^{h( \gamma,p)/2 )}}-
\frac{1}{ {\m p\m}^{h(\theta_0)}}\right|\m\rho(\gamma)\m
d\gamma\\
&+&\frac{1}{ {\m p \m}^{h( \theta_0 )}}\int_{\m \gamma \m \le \m p
\m^{s}}\m c(\gamma,p) -c(\theta_0)\m \m\rho(\gamma)\m d\gamma.
\end{eqnarray*}
\\
Let us use the Lipschitz assumptions on  $h$ and $c$.\\
\begin{Lem}\label{majorationlips}
If $g$ is an homogeneous function of degree $0$, Lipschitz of order $\alpha$ on the sphere
$S^{d-1}$,  then there exists $C>0$ such that for all $p\neq 0$ and $ \gamma \in \R^{d-1}$,
$$ \m g(\gamma,p)- g(0,p) \m \le C\min{\left(\left(\frac{\m \gamma \m}{\m p \m}\right)^{\alpha},1\right)}.$$
\end{Lem}
\begin{proof}
The function $g$ is continuous on the sphere and thus it is bounded. Then, for $p\neq 0$ and $ \gamma \in \R^{d-1}$,
$$\m g(\gamma,p)- g(0,p) \m \le 2\n g\n_{\infty}.$$
Moreover, $g$ is  Lipschitz of order $\alpha$ on the sphere. Then, there exists $C>0$ such that, for $p\neq 0$,
$$\m g(\gamma,p)- g(0,p) \m \le \left| \frac{(\gamma,p)}{(\m \gamma\m^2+p^2)^{1/2}}-\frac{(0,p) }{\m p\m}
 \right|^{\alpha}.$$
But,
$$
\left| \frac{(\gamma,p)}{(\m \gamma\m^2+p^2)^{1/2}}-\frac{(0,p) }{\m p\m}
 \right|^2 = 2\left(1-\left(1+\frac{\m \gamma\m^2}{p^2}\right)^{-1/2}\right)
\le \frac{\m \gamma\m^2}{p^2},$$
which concludes the proof of Lemma \ref{majorationlips}.
\end{proof}
\noindent
Let us recall that $c$ is an even homogeneous function thus $c(\theta_0)=c(0,p)$, for all $p\neq 0$.
Then, since $c$ is Lipschitz of order $\alpha$, one can find  $C_1>0$ such that
$$\Delta_{s,2}(p)=\frac{1}{{\m p\m}^{h(\theta_0)}}\int_{\m \gamma \m \le \m p
\m^{s}}\m c(\gamma,p)-c(0,p)\m \m\rho(\gamma)\m d\gamma
 \le C_1\m p\m^{-h(\theta_0)-\alpha (1-s)},$$
which is negligible compared to $\m p \m^{-h(\theta_0)-\delta}$ as soon as $\delta<\alpha(1-s)$.\\
\\
It remains to consider
\begin{eqnarray*}
\Delta_{s,1}(p)&=&\int_{\m \gamma \m \le \m p \m^{s}} \m c(\gamma,p)\m\left|
\frac{1}{ {\left(\m\gamma\m^2 +p^2\right)}^{h(\gamma,p)/2}}-
\frac{1}{ {\m p \m}^{h( \theta_0 )}}\right|\m \rho(\gamma)\m d\gamma\\
&=&\frac{1}{{\m p\m}^{h(\theta_0)}}\int_{\m \gamma \m \le \m p \m^{s}}
\m c(\gamma,p)\m\left|\frac{{\m p\m}^{h(\theta_0)}}{ {\left(\m\gamma\m^2 + p^2\right) }^{h(
\gamma,p)/2}}-1\right|\m \rho(\gamma)\m d\gamma.
\end{eqnarray*}
Let us write
$$\frac{{\m p\m}^{h(\theta_0)}}{ {\left(\m\gamma\m^2 +p^2\right)}^{h(
\gamma,p)/2}}=e^{l(p)}, $$ where, for $p\neq 0$,
\begin{eqnarray*}
l(p)&=&h(\theta_0)\ln{\m p\m}-\frac{1}{2}h(
\gamma,p)\ln{\left(\m p\m^2+\m\gamma\m^2\right)}\\
&=&\left(h(0,p)-h(\gamma,p)\right)\ln{\m p\m}-\frac{1}{2}h(
\gamma,p)\ln{\left(1+\frac{\m\gamma\m^2}{\m p\m^2}\right)},
\end{eqnarray*}
writing $h(\theta_0)=h(0,p)$ since $h$ is an even homogeneous function.
Since $h$ is Lipschitz of order $\alpha$, by Lemma \ref{majorationlips}, for $s<1$, there exists $C_2>0$ such that,
for $\m p\m\ge A>e$ and
$\m \gamma \m \le \m p \m^{s}$,
\begin{eqnarray*}
\m l(p)\m & \le &C_2\left(\left(\frac{\m \gamma \m}{\m p
\m}\right)^{\alpha}\ln{\m p\m}+\frac{\m\gamma\m^2}{\m p\m^2}\right)\\
&\le &2C_2\left(\frac{\m \gamma \m}{\m p
\m}\right)^{\alpha}\ln{\m p\m}
\le 2C_2{\m p\m}^{-\alpha (1-s)}\ln{\m p\m}.
\end{eqnarray*}
The function $t \mapsto {|t| }^{-\alpha(1-s)}\ln\m t \m$
 tends to $0$ at infinity. Thus one can find
$A_{s}>0$ such that for $\m p \m>A_{s}$ we get
$\m l(p)\m <1$. Then, for $\m p \m>A_{s}$
$$\m e^{l(p)} -1 \m \le e\m l(p)\m \le 2eC_2{\m p\m}^{-\alpha (1-s)}\ln{\m p\m},$$
and finally
$$\Delta_{s,1}(p)
\le 2eC_2\n c \n_{\infty}\m p \m^{-h(\theta_0)-\alpha (1-s)}\ln{\m p\m}.$$
For $\delta<\alpha$, since $\left|\Delta_s(p)\right|\le \Delta_{s,1}(p)+\Delta_{s,2}(p)$, taking $s\in(0,\frac{\alpha-\delta}{\alpha})\subset
(0,1)$ we get $$\Delta_s(p)=o_{\m p \m \rightarrow +\infty}\left(\frac{1}{\m p \m^{h(\theta_0)+\delta}}\right)$$
and (\ref{TRAFBc}) follows.
In the general case, let us assume
that, for all $\delta\in (0,\delta_0)$ and $\xi\in\R^d$,
$$f(\xi)= \frac{c(\xi)}{ {\m \xi \m}^{h(\xi)}}+ {o}\left(\frac{1}{\m \xi
\m^{h(\xi)+\delta}}\right) \mbox{ when }
{\m \xi \m \rightarrow +\infty}.$$
Replacing $\rho$ by $\m\rho\m$ and $h$ by $h+\delta$ in the special case above, we get the result for the remainder.
\end{proof}
Let us remark that to simplify the statement of this proposition,
 we assumed 
the window to be in
the Schwartz class. However it is proved in  \cite{HB} p. 85 that
 the result still holds for a window $\rho\in L^1(\R^{d-1})$ that satisfies
\begin{equation}\label{majoration de rho}
\left|\rho(\gamma)\right|={{O}}\left(\frac{1}{\m\gamma\m^{M+d-1}}\right)\mbox{ when } \gamma\in \R^{d-1}
\mbox{ and }{\m \gamma \m \rightarrow +\infty},
\end{equation}
with $M> h_1-h_0$ for $h\left(S^{d-1}\right)=[h_0,h_1]$. In this case $\delta_1=\min\left(\delta_0, \alpha\frac{M+h_0-h_1}{M +\alpha}\right)$.\\
\\

We can now state our main result concerning the identification of the anisotropic index of an anisotropic fractional Brownian field. We keep the notations of part \ref{GQV1D} for the generalized quadratic variations of a 1D-process and recall that we fix the direction $\theta_0=(0,\ldots,0,1)\in S^{d-1}$.
\begin{Th}\label{idAFB}
Let  $X=\{X(t);t\in \R^d\}$ be an anisotropic fractional Brownian field,
with anisotropic index $h$ given by an even homogeneous function of degree $0$ with values in $(0,1)$, which is assumed to be in ${\mathcal C}^1\left(S^{d-1}\right)$.\\
Let $\rho$ be a  window in ${\mathcal S}\left(\R^{d-1}\right)$.
Let $R_{\rho}X$ be the Radon transform of the field $X$  with the window $\rho$.\\
Let $\mathbf{a}$ be a filter of order $K$ and $u,v\ge 1$ two integers with $u\neq v$.  If  $K>h(\theta_0)+\frac{d-1}{2}$ then almost surely
$$\widehat{h_{N,\mathbf{a}}(\theta_0)}(u,v)=\frac{1}{2\log(u/v)}\log{\left(\frac{V_{N,\mathbf{a^u}}(R_{\rho}X)}{V_{N,\mathbf{a^v}}(R_{\rho}X)}\right)}-\frac{d-1}{2}\underset{N\rightarrow+\infty}\longrightarrow h(\theta_0),$$
Moreover, if $K>h(\theta_0)+\frac{d-1}{2}+1/4$
$$\sqrt{N}\left(\widehat{h_{N,\mathbf{a}}(\theta_0)}(u,v)-h(\theta_0)\right)\overset{d}{\underset{N\rightarrow +\infty}{\longrightarrow}}{\mathcal N}\left(0,\gamma_\mathbf{a}^{u,v}\left(h(\theta_0)+\frac{d-1}{2}\right)
\right),$$
with
$$N\mathbb{E}\left(\left(\widehat{h_{N,\mathbf{a}}(\theta_0)}(u,v)-h(\theta_0)\right)^2\right)\underset{N\rightarrow +\infty}{\longrightarrow} \gamma_\mathbf{a}^{u,v}\left(h(\theta_0)+\frac{d-1}{2}\right).$$
\end{Th}
\begin{proof}
 It is sufficient to prove that the spectral density of $R_{\rho}X$ satisfies the assumption of Proposition \ref{propVNDAS}. From (\ref{TRDS}),
this spectral density is given by the function
 $${\mathcal R}_{\left|\widehat{\rho}\right|^2}f_h(p)=\int_{\R^{d-1}}f_h((\gamma,p))\left|\widehat{\rho}(\gamma)\right|^2d\gamma,\mbox{ for all } p\in \R,$$
with $f_h$ given by (\ref{density}).
Since we can divide $R_{\rho}X$ by a constant, we can assume that $\int_{\R^{d-1}}
\left|\widehat{\rho}(\gamma)\right|^2d\gamma=1$. Then, since $2h+d$ satisfies assumptions of Proposition \ref{DASRT} with $\alpha=1$, 
by Proposition \ref{DASRT}, for all $\delta<1$,
 $${\mathcal R}_{\left|\widehat{\rho}\right|^2}f_h(p) =\frac{1}{ {\m  p \m}^{2h(\theta_0)+d}}+
\underset{\m p \m \rightarrow +\infty}{o}\left(\frac{1}{\m p \m^{2h(\theta_0)+d+\delta}}\right).$$
We can also write this as
$${\mathcal R}_{\left|\widehat{\rho}\right|^2}f_h(p) =\frac{1}{\m p \m^{2\left(h(\theta_0)+\frac{d-1}{2}\right)+1}}+
\underset{\m p \m \rightarrow +\infty}{o}\left(\frac{1}{\m p \m^{2\left(h(\theta_0)+\frac{d-1}{2}\right)+1+\delta}}\right).$$
Moreover, since $h\in {\mathcal C}^1\left(S^{d-1}\right)$, the function
$g(\gamma):=h(\gamma,1)$ is differentiable on $\R^{d-1}$ and, for $\gamma, x\in\R^{d-1}$,
$$D_{\gamma}g(x)=\frac{1}{(\m\gamma\m^2+1)^{1/2}}D_{\frac{(\gamma,1)}{\m(\gamma,1)\m}}h\left((x,0)-\frac{\gamma
\cdot x
}{\m\gamma\m^2+1}(\gamma,1)\right),$$
such that $$\left|D_{\gamma}g(x)\right|\le 2\n Dh\n_{\infty}\m x\m.$$
Thus,
the spectral density $f_h$ is differentiable on $\R^{d}\smallsetminus\{0\}$. Let $(\gamma,p)\in \R^{d-1}\times\R^*$, since $h(\gamma,p)=g(\gamma/p)$, we get
$$\frac{\partial}{\partial p}f_h(\gamma,p)=f_h(\gamma,p)\left(- D_{\frac{\gamma}{p}}g\left(-\frac{\gamma}{p^2}\right)\ln\left(\m\gamma\m^2+p^2\right)- \frac{p\left(2h(\gamma,p)+d\right)}{\m\gamma\m^2+p^2}\right).$$
It follows that the spectral density ${\mathcal R}_{\left|\widehat{\rho}\right|^2}f_h$ is differentiable on $\R\smallsetminus\{0\}$ and for $p\neq 0$ we have
$$\left({\mathcal R}_{\left|\widehat{\rho}\right|^2}f_h\right)'(p)=\int_{\R^{d-1}}\frac{\partial}{\partial p}f_h(\gamma,p)\left|\widehat{\rho}(\gamma)\right|^2d\gamma.$$
Let us write
$$\frac{\partial}{\partial p}f_h(\gamma,p)=\frac{\m \gamma\m}{p^2}F_1(\gamma,p)-pF_2(\gamma,p),$$
with 
$$F_1(\gamma,p)=D_{\frac{\gamma}{p}}g\left(-\frac{\gamma}{|\gamma|}\right)\ln\left(\m\gamma\m^2+p^2\right)f_h(\gamma,p)$$
and
$$F_2(\gamma,p)= \left(2h((\gamma,p))+d\right)f_{h+1}(\gamma,p).$$
Therefore,
$$\left({\mathcal R}_{\left|\widehat{\rho}\right|^2}f_h\right)'(p)=\frac{1}{p^2}{\mathcal R}_{\m\gamma\m\left|\widehat{\rho}\right|^2}F_1(p)-p{\mathcal R}_{\left|\widehat{\rho}\right|^2}F_2(p).$$
Then, by Proposition \ref{DASRT}, whenever $\delta<1$,
$${\mathcal R}_{\left|\widehat{\rho}\right|^2}F_2(p)=
\frac{2h(\theta_0)+d}{ {\m  p \m}^{2h(\theta_0)+d+2}}+
\underset{\m p \m \rightarrow +\infty}{o}\left(\frac{1}{\m p \m^{2h(\theta_0)+d+2+\delta}}\right).$$
Moreover, for any $\epsilon>0$ small enough
$$\left|F_1(\gamma,p)\right|\le 2\n Dh\n_{\infty}f_{h-\epsilon}(\gamma,p).$$
Since $|\gamma|\left|\widehat{\rho}\right|^2$ is integrable over $\R^{d-1}$ and rapidly decreasing, following the same lines as in the proof of Proposition \ref{DASRT}, we get
$$\left|{\mathcal R}_{\m\gamma\m\left|\widehat{\rho}\right|^2}F_1(p)\right|\le C{\m  p \m}^{-2h({\theta_0})-d+2\epsilon}.$$
This allows us to conclude that
$$\left({\mathcal R}_{\left|\widehat{\rho}\right|^2}f_h\right)'(p)=-
\frac{2h(\theta_0)+d}{ {\m  p \m}^{2h(\theta_0)+d+1}}+\underset{\m p \m \rightarrow +\infty}{o}\left(\frac{1}{ {\m  p \m}^{2h(\theta_0)+d+1}}\right).$$
Therefore $R_{\rho}X$ satisfies the assumption of Proposition \ref{propVNDAS} with
$H=h(\theta_0)+\frac{d-1}{2}$ and $s<1$, which concludes the proof.
\end{proof}
Let us remark that if we choose a window such that
$\left|\widehat{\rho}\right|^2$ only satisfies (\ref{majoration de rho}) for $M>2(h_1-h_0)$ with $h\left(S^{d-1}\right)=[h_0,h_1]$,
the estimator $\widehat{h_{N,\mathbf{a}}(\theta_0)}(u,v)$ still tends almost surely to $h(\theta_0)$.
However the speed of convergence will depend on $M$. Actually the result of Theorem \ref{idAFB}
holds when  $M\ge 1+4(h_1-h_0)$ whereas for $M< 1+4(h_1-h_0)$, Proposition \ref{propVNDAS} shows that, for all $s<\frac{M-2(h_1-h_0)}{M +1}$,
$$N^{s}\left(\widehat{h_{N,\mathbf{a}}(\theta_0)}(u,v)-h(\theta_0)\right)\overset{L^2(\Omega)}{\underset{N\rightarrow +\infty}{\longrightarrow}}0.$$
\\

Let us point out that we have restricted
$h$ to have values in $(0,1)$ so that the anisotropic fractional Brownian field is well defined. This can be weakened by considering fields with higher order stationary increments or with spectral density asymptotically of the order of $f_h$. Similar results can be obtained with the further assumption that
the partial derivatives of order $1$ are
asymptotically of the order of the partial derivatives of  $f_h$
(see \cite{HB} for example).\\

For numerical applications, one has to approximate the Radon transform of $X$. Following \cite{ABAEAA}
we replace $R_{\rho}X(t)$ by
\begin{equation}\label{discretisation}
I_M(t)=M^{-(d-1)}\sum_{s\in \Z^{d-1}}
 X\left(M^{-1}s,t\right) \rho(M^{-1}s),
 \end{equation}
for $\rho$ a smooth window with compact support and $M\ge 1$ an integer. Let us denote for $\mathbf{a}$ a filter of order $K$ and length $l$ and $N\le M$,
\begin{equation} \label{projestimator}
T_{N,\mathbf{a}}^M(X)=\frac{1}{N-l+1}\sum_{p=0}^{N-l} \left( \sum_{k=0}^{l}a_kI_M\left(\frac{p+k}{N}\right)
\right)^2.
\end{equation}
The key point of the proof is to estimate the error due to the approximation of $V_{N,\mathbf{a}}\left(R_{\rho}X\right)$ by $T_{N,\mathbf{a}}^M(X)$. Let us denote $H=h(\theta_0)+\frac{d-1}{2}$. Under the assumptions of Theorem \ref{idAFB},
following the same lines as in \cite{ABAEAA}, for $\alpha>0$ with $\alpha<h_0=\underset{S^{d-1}}{\min h}$, one can prove that there exists a positive finite random variable $C$ such that, for all $N\ge 1$,
$$\left|T_{N,\mathbf{a}}^M(X)^{1/2}-V_{N,\mathbf{a}}\left(R_{\rho}X\right)^{1/2}\right|\le CM^{- \alpha} \mbox{ a.s.}$$
Moreover, since for $K>H$ by Proposition \ref{esperance}  $\frac{V_{N,\mathbf{a}}\left(R_{\rho}X\right)}{\mathbb{E}\left(V_{N,\mathbf{a}}\left(R_{\rho}X\right)\right)}\underset{N\rightarrow +\infty}{\longrightarrow}1$ a.s., 
one can find a positive finite random variable $C'$ such that
$$V_{N,\mathbf{a}}\left(R_{\rho}X\right)^{-1/2}\le C'N^{H} \mbox{ a.s.}$$
Then, for $M^{-\alpha}N^H$ small enough, writing
$$\log\left(\frac{T_{N,\mathbf{a}}^M(X)}{V_{N,\mathbf{a}}\left(R_{\rho}X\right)}\right)=2\log\left(1+\frac{
T_{N,\mathbf{a}}^M(X)^{1/2}-V_{N,\mathbf{a}}\left(R_{\rho}X\right)^{1/2}}{V_{N,\mathbf{a}}\left(R_{\rho}X\right)^{1/2}}\right),$$
one can find a positive finite random variable $C''$ such that, a.s.
$$\left|\log\left(\frac{T_{N,\mathbf{a}}^M(X)}{V_{N,\mathbf{a}}\left(R_{\rho}X\right)}\right)\right|\le C''M^{-\alpha}N^H.$$
We can state the following result.
\begin{Prop}\label{idstep} We keep the assumptions of Theorem \ref{idAFB} and take $\rho$ with compact support. Let  $H=h(\theta_0)+\frac{d-1}{2}$ and $h_0=\underset{S^{d-1}}{\min h}$.
 If  $K>H$ and $q>H/h_0$ then, almost surely,
$$\widehat{H_{N,\mathbf{a}}(\theta_0)}(u,v)=\frac{1}{2\log(u/v)}
\log{\left(\frac{T_{N,\mathbf{a^u}}^{N^q}(X)}{T_{N,\mathbf{a^v}}^{N^q}(X)}\right)}-\frac{d-1}{2}\underset{N\rightarrow+\infty}\longrightarrow h(\theta_0),$$
Moreover,  if  $K>H+1/4$ and $q>(H+1/2)/h_0$  then
$$\sqrt{N}\left(\widehat{H_{N,\mathbf{a}}(\theta_0)}(u,v)-h(\theta_0)\right)\overset{d}{\underset{N\rightarrow +\infty}{\longrightarrow}}{\mathcal N}\left(0,\gamma_\mathbf{a}^{u,v}(H)\right).$$
\end{Prop}

\section{Numerical Study}

In this section we present a preliminary evaluation of projection-based estimators studied above. We first describe the synthetic datasets used for our evaluation and then discuss our estimation results.  

\subsection{Simulation}

A lot of numerical methods have been proposed these last years to simulate 1-dimensional fractional Brownian motion (fBm). Most of them give rise to approximate syntheses, such as the midpoint displacement method (see \cite{Norros}, for instance), the wavelet based decomposition (\cite{MST}, \cite{Abry}, \cite{Pipiras}, etc...), or more recently a method based on correlated random walks \cite{Enriquez}. A few of them can be applied not only for 1-dimensional fBm but also to simulate 2-dimensional (anisotropic) fractional Brownian fields and  lead to approximate syntheses. Of course there exist exact synthesis methods based on the Choleski decomposition of the covariance function. This yields numerical problems due to the size of the matrix. In order to have fast synthesis one can use the stationarity of the increments   by applying the embedding circulant matrix method \cite{Dietrich}. By this way, we easily obtain fast and exact synthesis of 1-dimensional fBm \cite{Ileana}. Some authors, as in \cite{Chan} and \cite{KaplanKuo}, apply this method for higher dimension but this does not yield to exact synthesis. Finally, M. L. Stein proposed a fast and exact synthesis method for isotropic fBm surfaces in \cite{Stein}.

For a preliminary evaluation, we used this method for generating a dataset containing three subsets of $1000$ fBm fields simulated with each of the Hurst parameter values: $h=0.2$ (low regularity) $h=0.5$ (medium regularity) $h=0.7$ (high regularity). fBm fields were generated on a discrete grid $G$ of size $(M+1) \times (M+1)$ ($M=2^9=512$) which is defined by
\begin{equation} \label{egn:grid}
G = \bigg\{ \bigg(\frac{k_1}{M},\frac{k_2}{M} \bigg), 0 \leq  k_1,k_2 \leq M\bigg\}. 
\end{equation}
The method was implemented in collaboration with A. Fraysse and C. Lacaux, and the corresponding matlab codes are available at  \url{http://ciel.ccsd.cnrs.fr}. 

%With this algorithm we generate realizations of points of fBm of Hurst parameter $H \in (0,0.7)$ on an equispaced grid with mesh $M_0=2^{10}$, namely $$\left\{B_H\left(\frac{k}{M_0}, \frac{l}{M_0}\right); 0\le k,l\le M_0\right\},$$ where $M_0=764$ is the greatest integer less than $M_0/\sqrt{2}$.\\ In order to perform a Radon transform, we have thus chosen $\rho=\mathbf{1}_{[0,M_0/M_0]}$. Let us remark that the Fourier transform of this window is given, for $\gamma\in\R$, by $$\widehat{\rho}(\gamma)=\frac{e^{-iM_0\gamma/M_0}-1}{-i\gamma},$$ such that $\left|\widehat{\rho}\right|^2$ satisfies (\ref{majoration de rho}) for any $M<1$.

\begin{figure}[th]
\begin{center}
\begin{tabular}{ccc}
\includegraphics[width=0.24\textwidth, angle=-90]{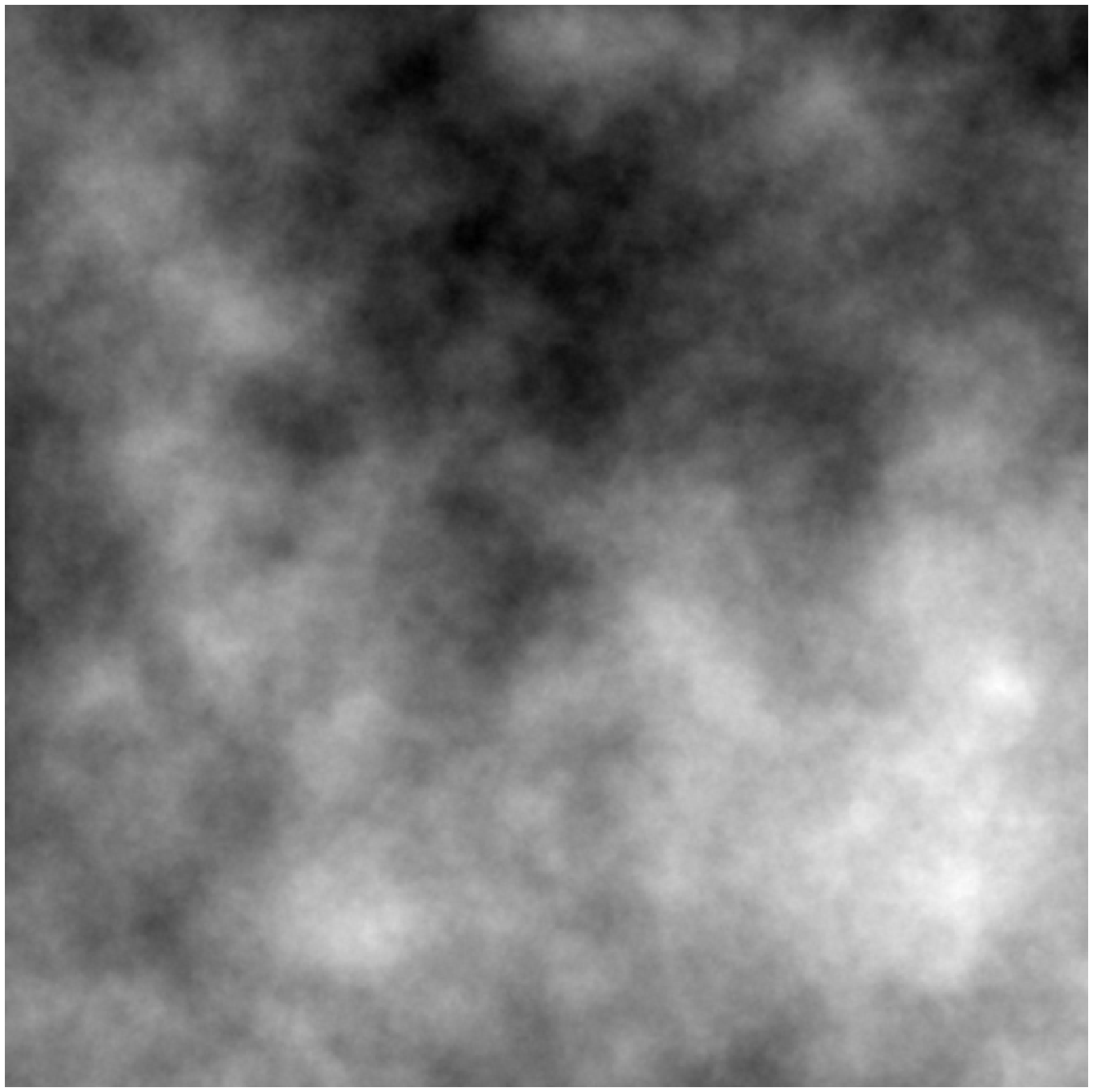} & \includegraphics[width=0.24\textwidth, angle=-90]{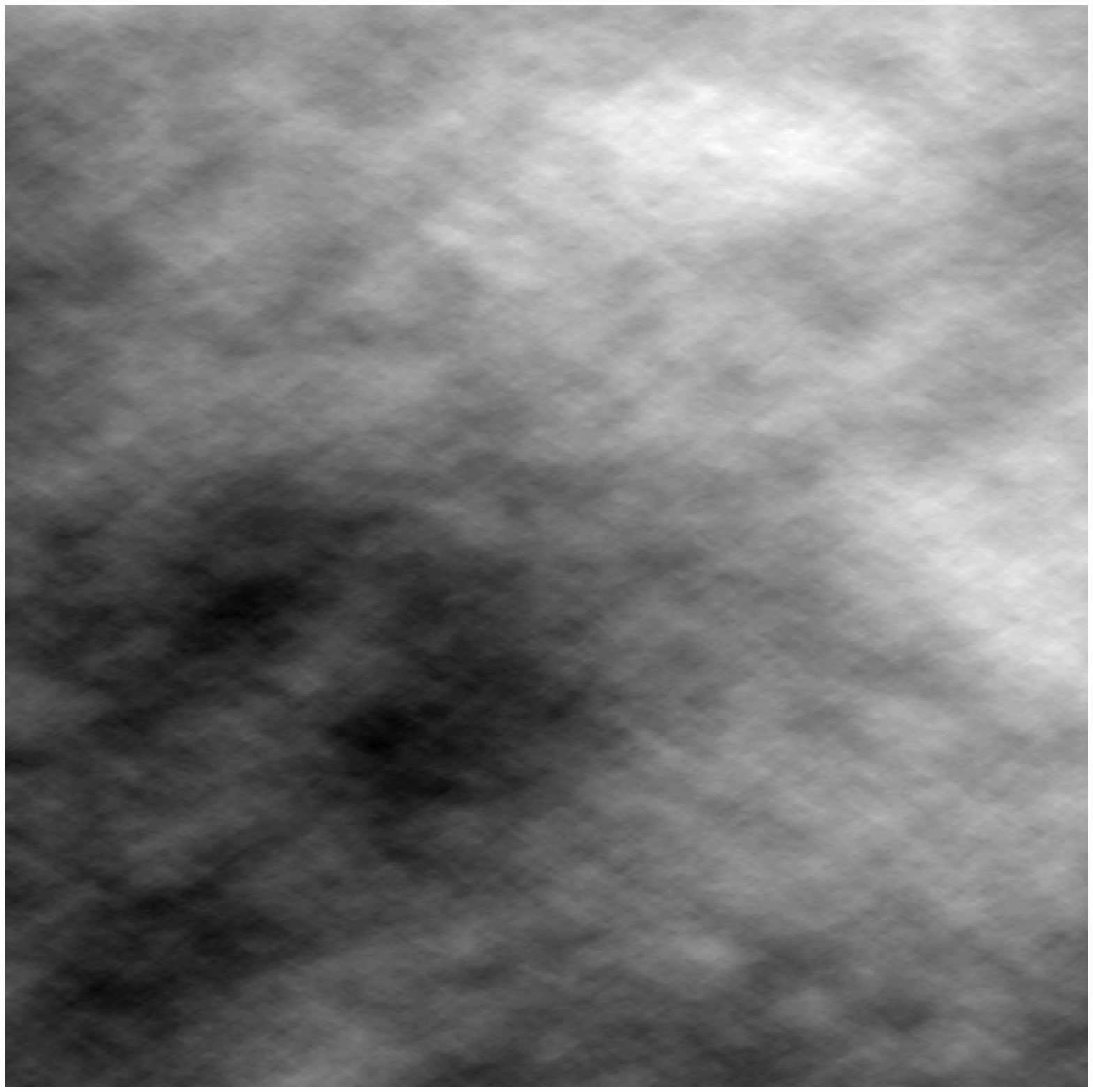} & \includegraphics[width=0.24\textwidth, angle=-90]{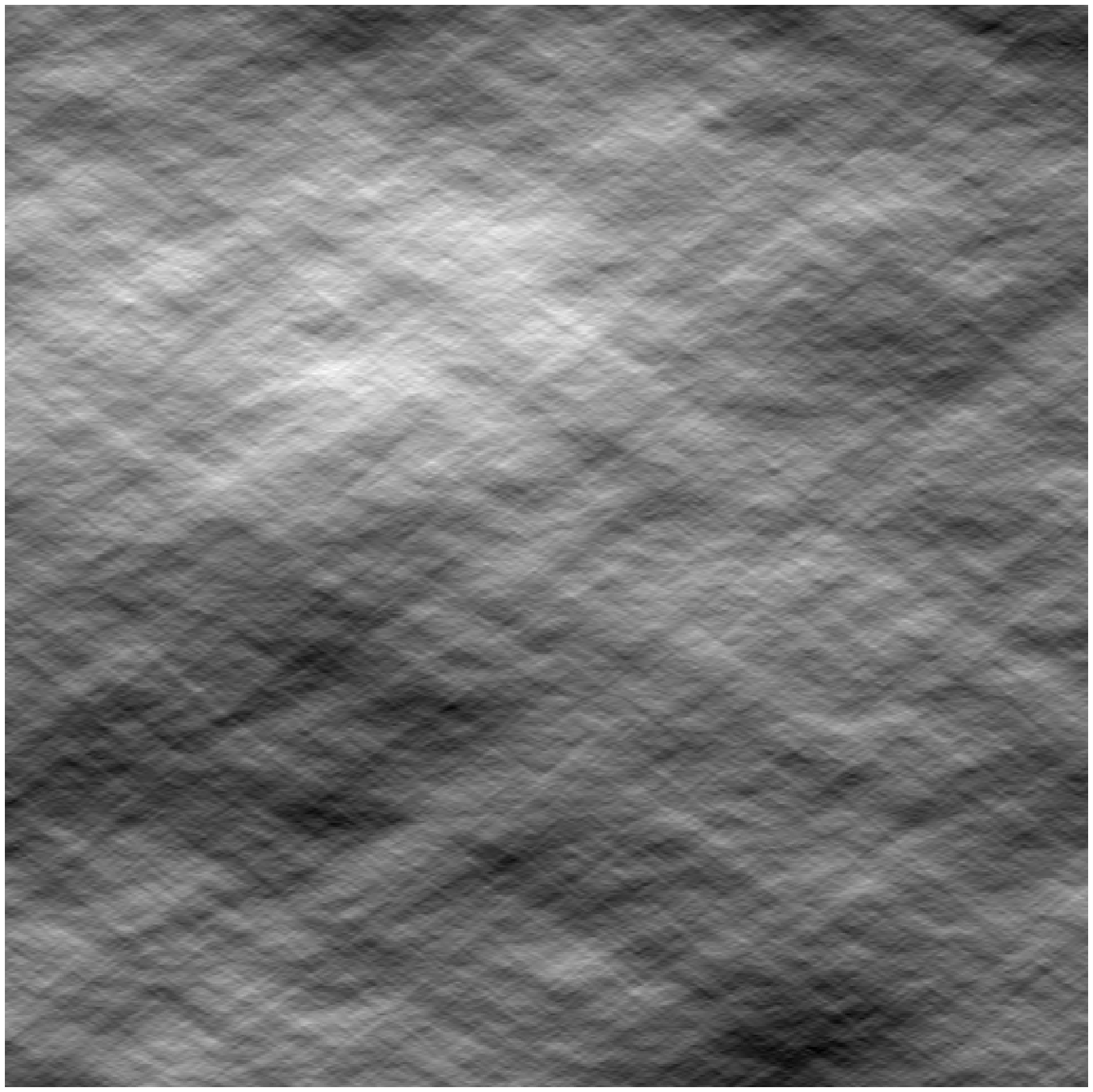} \\
(a) & (b) & (c)
\end{tabular}
\end{center} \caption{Some simulations of afB surfaces. Horizontal and vertical regularity parameters used for simulations are  $0.7$ and $0.7$ in (a), $0.7$ and $0.5$ in (b) and $0.7$ and $0.2$ in (c).}\label{examples}
\end{figure}

However, this method can not be generalized to simulate anisotropic fields. Hence we also used a spectral representation approximation (SRA) technique to generate anisotropic fractional Brownian fields (afB). The regularity of two-dimensional  afB we generated differs in both vertical and horizontal directions. The spectral densities of these fields are of the form
\begin{equation}
\forall  \: \xi=(\xi_1,\xi_2) \in \mathbb{R}^2\smallsetminus\{(0,0)\}, \: f(\xi) = \left\{
\begin{array}{ll} 
\vert \xi \vert^{-(2 h_v +2)}, & \mathrm{if} \:\: \vert \xi_1 \vert < \vert \xi_2 \vert \\
\vert \xi \vert^{-(2 h_h +2)}, & \mathrm{otherwise}.\\
 \end{array}  
\right.
\end{equation}
In this expression, $h_h$ and $h_v$ form a pair of parameters in $(0,1)$ which characterize the anisotropy of  generated fields. Their regularity is then given by $\min(h_v,h_h)$. We used an approximation of the spectral representation in Equation (\ref{spectral})  which was obtained by discretization. 
The discrete approximation is given by 
$$x\left(\frac{k_1}{M},\frac{k_2}{M}\right)=\Re\left(y\left(\frac{k_1}{M},\frac{k_2}{M}\right)-y(0,0)\right), \mbox{ for } 0\le k_1, k_2\le M,$$
with
\begin{equation} \label{eqn:spectral_discrete} 
y\bigg(\frac{k_1}{M},\frac{k_2}{M}\bigg)= \pi\sum_{n_1= - M+1}^{M} \sum_{n_2 = - M+1}^{M}  z(n_1,n_2) \: g\bigg(\pi n_1, \pi n_2\bigg)  e^{-\frac{ 2i\pi }{2M} (n_1 k_1+n_2 k_2)}, 
\end{equation} 
where $z(n_1,n_2)$ are $ (2M)^2$ independant realizations of complex random variables whose real and imaginary components are two uncorrelated zero-mean standard gaussian variables and the function $g(x,y)=f^{1/2}(x,y)$ for all $(x,y) \in \R^2\smallsetminus\{(0,0)\}$ and $g(0,0)=0$. %Setting $T^2=8 \pi M^2$ and letting $M$ tends to infinity, the limit of the sum in Equation (\ref{eqn:spectral_discrete}) is a riemman series converging to the integral defined in Equation (\ref{spectral}). 
From a practical point of view, the sum in Equation  (\ref{eqn:spectral_discrete}) can be interpreted as a filtering in the Fourier domain of a white noise $z$ by a low-pass linear filter characterized by the transfer function $g$. Based on this method, afB fields approximations can be easily and quickly simulated using the fast Fourier transform. Similar simulations are used in \cite{LESI} for an evaluation of a different estimator.  Such an approach can also be  extended to the simulation of 3D fields and anisotropic fields having more than two different directional regularities. As an illustration, some simulation examples are shown on Figure \ref{examples}.

By this way, we generated a second evaluation dataset containing six subsets of  $1000$ fields on the grid $G$ for each of the following regularity parameter pairs: low isotropic regularity $(h_h=0.2,h_v=0.2)$, medium isotropic regularity $(h_h=0.5,h_v=0.5)$, high isotropic regularity $(h_h=0.7,h_v=0.7)$, high/medium anisotropic regularity  $(h_h=0.7,h_v=0.5)$, high/low anisotropic regularity $(h_h=0.7,h_v=0.2)$, medium/low anisotropic regularity $(h_h=0.5,h_v=0.2)$.

\subsection{Estimation}

Given a  field realization $x$, we first computed discrete Radon transforms $I_h$  and $I_v$ in both horizontal and vertical directions, as defined in Equation (\ref{discretisation}) for a window function of the shape $\mathbf{1}_{[0,1]}$: for all $0\le k ,l \le M$,
\begin{equation}
 I_h \bigg(\frac{k}{M}\bigg) = \frac{1}{M} \sum_{k_2=0}^{M} x\bigg(\frac{k}{M},\frac{k_2}{M}\bigg) \:\: \mathrm{and} \:\:  I_v  \bigg(\frac{l}{M}\bigg) =  \frac{1}{M} \sum_{k_1=0}^{M} x\bigg( \frac{k_1}{M},\frac{l}{M} \bigg).
\end{equation}
%Let us remark that the Fourier transform of this window is given, for $\gamma\in\R$, by $$\widehat{\rho}(\gamma)=\frac{e^{-i \gamma}-1}{-i\gamma}$$, so that $\left|\widehat{\rho}\right|^2$ satisfies (\ref{majoration de rho}) for any $M<1$.

Then we computed generalized quadratic variations $T^\nu_{e,1}$ for both directions ($e=v,h$), as defined in Equation (\ref{projestimator}), with the second order filter  $a^1=(1,-2,1)$  and a step of size $1/N$ where $N=M/2^\nu$,  for $\nu$ in $\{0, 1, 2, 3\}$:
\begin{equation}
T^\nu_{e,1}=\frac{1}{(M/2^{\nu}-1)}\sum_{p=0}^{M/2^\nu-2} \bigg( I_e \bigg( 2^\nu \frac{p}{M} \bigg) - 2 \: I_e \bigg( 2^\nu \frac{p+1}{M} \bigg) +  I_e \bigg( 2^\nu \frac{p+2}{M}\bigg)  \bigg)^2,
\end{equation} 
We also computed generalized quadratic variations $T^\nu_{e,2}$ with the dilated second order filter $a^2=(1,0,-2,0,1)$: 
\begin{equation}
T^\nu_{e,2}=\frac{1}{(M/2^{\nu}-3)}\sum_{p=0}^{M/2^{\nu}-4}  \bigg( I_e\bigg( 2^\nu \frac{p}{M} \bigg) - 2 \: I_e \bigg( 2^\nu \frac{p+2}{M} \bigg) +  I_e \bigg( 2^\nu \frac{p+4}{M} \bigg)   \bigg)^2.
\end{equation} 

Finally, we obtained the projection-based estimate $\hat{h}^\nu_e$ of the  index $h_e$ of $X$ in the  direction  $e$ ($e=v,h$) for $\nu \in \{0,\cdots,3\}$ as
\begin{equation}
\hat{h}^\nu_e = \frac{1}{2 \log(2)} \log \bigg( \frac{T^\nu_{e,2}}{T^\nu_{e,1} }\bigg)-\frac{1}{2}.
\end{equation}

For the evaluation, we computed empirical biases and standard deviations of estimators over sets of simulated fields having same characteristics. Empirical biases were obtained as a difference between the mean parameter estimates and the real parameter value. In order to enhance the estimator ability to capture anisotropic properties, we also computed biases and standard deviations of differences between horizontal and vertical estimates. Whole results are reported in Tables \ref{results} and \ref{results2}.

\begin{table}[th] 
\begin{center}
\begin{tabular}{|c||c||r@{$\pm$}l|r@{$\pm$}l|r@{$\pm$}l|}
\hline
$h$ &  $\nu$ & $b_h$ &  $\sigma_h$ & $b_v$ & $\sigma_v$ & $b_{h,v}$ &  $\sigma_{h,v}$ \\
\hline\hline 
  0.7 &  0 & -0.047 &  0.049  &-0.045 &  0.045  &-0.002 &  0.069 \\
  0.7 &  1 &-0.012  &  0.056  &-0.017 &  0.061  &  0.005  &  0.085 \\
  0.7 &  2 &-0.03   &   0.081 &-0.023 &  0.093  &-0.007 &  0.124 \\
  0.7 &  3 &-0.054  &  0.113  &-0.04  &   0.114 &-0.014 &  0.158 \\
\hline
 0.5 &  0 &-0.092 &  0.052 &-0.095 &  0.052 &  0.003 &  0.073\\
 0.5 &  1 &-0.034 &  0.069 &-0.035 &  0.071 &  0.001 &  0.099\\
 0.5 &  2 &-0.007 &  0.099 &-0.031 &  0.095 &  0.024 &  0.132\\
 0.5 &  3 &-0.029 &  0.131 &-0.023 &  0.136 &-0.006 &  0.192\\
\hline 
0.2 &  0 & -0.239 &  0.055  &-0.245 &  0.062 &  0.006 &  0.082\\
  0.2 &  1&  -0.112 &  0.079 & -0.131 &  0.082 &  0.019 &  0.119\\
   0.2 &  2 & -0.041 &  0.113 & -0.039 &  0.12 & -0.002 &  0.162\\
   0.2 &  3 &  0.002 &  0.163  &-0.039 &  0.147 &  0.041 &  0.238\\
\hline
\end{tabular}
\end{center}
\caption{Evaluation of directional regularity estimators on synthetic fBm surfaces of Hurst index $h$ simulated using the Stein's method of exact synthesis. Values  $b_h$ and $b_v$ are empirical biases of horizontal and vertical regularity estimators and  values $\sigma_h$ and $\sigma_v$ are their associated standard deviations. Values $b_{h,v}=b_h-b_v$ are differences between bias. Values $\sigma_{h,v}$ are standard deviations of differences between horizontal and vertical regularity estimates. }\label{results}
\end{table}

On Table \ref{results}, we observe that standard deviations of estimation errors increase as the subsampling factor $2^\nu$ increase, meaning that subsambling of the projected signal has an effect on the estimator stability. However,  standard deviations do not vary significantly when the subsampling factor is fixed and parameter values are changed. For instance, when $\nu=0$, standard deviations are $0.049$ and $0.055$ when parameters values are $0.7$ and $0.2$, respectively. The order of standard deviation variations is about $10^{-2}$, for any fixed subsampling factor.

Estimators underestimate the real parameter value. The underestimation increases as the parameter value decreases. For instance, estimation biases obtained for  $\nu=0$ are $-0.047$, $-0.092$ and $-0.239$ when  parameter values are $0.7$, $0.5$ and $0.2$, respectively. Subsampling the projected signal  reduces the underestimation bias. For instance, when $h=0.2$, the bias is reduced to $-0.041$ when the subsampling factor is increased to $2^2$.

\begin{table}[bh] 
\begin{center}
\begin{tabular}{|c|c||c||r@{$\pm$}l|r@{$\pm$}l|r@{$\pm$}l|}
\hline
$h_h$ & $h_v$ &  $\nu$ & $b_h$ &  $\sigma_h$ & $b_v$ & $\sigma_v$ & $b_{h,v}$ &  $\sigma_{h,v}$ \\
\hline\hline 
    0.7 &   0.7  &  0   & 0.068  &  0.041  &  0.069 &   0.041  &-0.001   & 0.06\\
    0.7 &   0.7  &  1 &   0.003  &  0.063&    0.   &    0.059 &   0.003  &  0.087\\
    0.7 &   0.7  &  2 &-0.012  &  0.09 & -0.014 &   0.087 &   0.002  &  0.126\\
    0.7  &  0.7 &   3 &-0.021  &  0.125& -0.024 &   0.13  &   0.003  &  0.182\\
\hline
    0.5  &  0.5  &     0 &   0.1   &   0.046 &   0.102 &   0.044 &-0.002 &   0.065\\
    0.5 &   0.5 &   1 &   0.012 &   0.07  &   0.009 &   0.067 &   0.003 &   0.097\\
    0.5 &   0.5 &   2 &-0.007 &   0.1   &-0.008 &   0.097 &   0.001 &   0.139\\
    0.5 &   0.5 &   3 &-0.013 &   0.142 &-0.015 &   0.146 &   0.002 &   0.207\\
\hline
    0.2  &  0.2  &  0&   0.156  &  0.052  &  0.16    & 0.05   &-0.004 &   0.073 \\
    0.2  &  0.2 &   1&   0.034  &  0.08   &  0.03    & 0.078  &  0.004 &   0.112\\
    0.2  &  0.2 &   2&    0.004 &   0.113 &   0.004  &  0.112 &   0.   &    0.158\\
    0.2  &  0.2 &   3& -0.004 &   0.163 &-0.007  &  0.164 &   0.003&    0.238\\
\hline
    0.7 & 0.5 &  0 &  0.071 &  0.041 &  0.1 &   0.044  &-0.029 &  0.061\\
   0.7 &  0.5 &  1 &  0.001 &  0.064 &  0.002 &  0.067  &-0.001 &  0.095\\
   0.7 &  0.5 &  2 & -0.015 &  0.089&  -0.005 &  0.101&  -0.01 &   0.133\\
   0.7 &  0.5 &  3 & -0.026 &  0.131&  -0.014 &  0.137&  -0.012 &  0.189\\
 \hline  
0.7 &  0.2 &  0 &  0.072 &  0.041 &  0.157 &  0.052 & -0.085 &  0.065\\
   0.7 &  0.2 &  1 & -0.002 &  0.061 &  0.029 &  0.078 & -0.031 &  0.1\\
   0.7 &  0.2 &  2 & -0.014 &  0.087 &  0.01 &   0.114 & -0.024 &  0.14\\
   0.7 &  0.2 &  3 & -0.022 &  0.128 & -0.009 &  0.163&  -0.013 &  0.21\\
\hline
  0.5 &  0.2 &  0 &  0.098 &  0.045 &  0.159 &  0.053 & -0.061 &  0.069\\
   0.5 &  0.2 &  1 &  0.006 &  0.072 &  0.032 &  0.079&  -0.026 &  0.108\\
   0.5 &  0.2 &  2 & -0.003 &  0.103 &  0.007 &  0.117&  -0.01 &   0.16\\
   0.5 &  0.2 &  3 & -0.002 &  0.142 & -0.009 &  0.163 &  0.007 &  0.211\\
\hline
\end{tabular}
\end{center}
\caption{Evaluation of directional regularity estimators on synthetic afB surfaces simulated using the SRA  method of approximate synthesis. Values $h_h$ and $h_v$ are horizontal and vertical regularity parameters used for field simulations. Values  $b_h$ and $b_v$ are empirical biases of horizontal and vertical regularity estimators and  values $\sigma_h$ and $\sigma_v$ are their associated standard deviations. Values $b_{h,v}=b_h-b_v$ are differences between bias. Values $\sigma_{h,v}$ are standard deviations of differences between horizontal and vertical regularity estimates. }\label{results2}
\end{table}

On Table \ref{results2}, we observe that estimation biases obtained on fields simulated using the SRA technique differ from those of Table \ref{results}. Indeed, these biases indicate that estimators overestimate parameter values. This is probably due to the SRA technique which generates fields which are smoother than what they should be. The best evaluation of estimator biases are those of Table \ref{results} which were obtained on exact simulations of fBm. However, we can notice that standard deviations of Table \ref{results2} are consistent with those of Table \ref{results}. Despite simulation errors of the SRA technique, we can rely on results of Table  \ref{results2} to get a sense of the estimator properties when fields are anisotropic. In particular, we see that estimator biases and standard deviations do not vary significantly from isotropic to anisotropic simulations. When the subsampling factor is fixed, standard deviations are about the same for isotropic and anisotropic simulations.  On anisotropic simulations, biases still vary when estimated parameter values are changed. But, the bias of a regularity estimate in one direction does not depend on the parameter value estimated in the other direction. For instance, when $\nu=0$, biases for the estimation of parameter value $h=0.2$ are about $0.15$ in all the simulation cases involving an index value of $0.2$ ($(h_h,h_v)=(0.2,0.2),(0.7,0.2),(0.5,0.2)$). In conclusion, the field anisotropy seems not to have any effects on the estimator stability.

Besides, we observe  on the last column of Tables \ref{results} and  \ref{results2} that standard deviations of estimate differences are about the same when $\nu$ is fixed. In particular, when $\nu=0$, standard deviations are about $0.065$. On isotropic cases, biases are about $0$. This suggests that it is possible to distinguish between isotropic fields and anisotropic fields for which absolute differences between horizontal and vertical regularities are above $0.065$.

\section*{Acknowledgements}
The authors  would  like  to warmly  thank Anne Estrade   for
 her relevant contribution  as well as for
very fruitful discussions, Aline  Bonami for many remarks and comments  simplifying  many computations lines.
They are also very grateful to the referees for improving the first version of this text.

\addcontentsline{toc}{chapter}{Bibliographie}
\bibliographystyle{plain}
\bibliography{Biblio}

\begin{thebibliography}{10}

\bibitem{Abry}
P.~Abry and F.~Sellan.
\newblock The wavelet-based synthesis for fractional {B}rownian motion proposed
  by {F}. {S}ellan and {Y}. {M}eyer: remarks and fast implementation.
\newblock {\em Appl. Comput. Harmon. Anal.}, 3:377--383, 1996.

\bibitem{ABAEAA}
A.~Ayache, A.~Bonami, and A.~Estrade.
\newblock Identification and series decomposition of anisotropic {G}aussian
  fields.
\newblock {\em {\sl Proceedings of the Catania ISAAC05 congress}}, 2005.

\bibitem{BardetE}
J.~M. Bardet, G.~Lang, G.~Oppenheim, A.~Philippe, S.~Stoev, and M.~S. Taqqu.
\newblock Semi-parametric estimation of the long-range dependence parameter: a
  survey.
\newblock In {\em Theory and applications of long-range dependence}, pages
  557--577. Birkhäuser Boston, 2003.

\bibitem{benassiFWN}
A.~Benassi, S.~Cohen, J.~Istas, and S.~Jaffard.
\newblock Identification of filtered white noises.
\newblock {\em Stochastic Process. Appl.}, 75(1):31--49, 1998.

\bibitem{BJR}
A.~Benassi, S.~Jaffard, and D.~Roux.
\newblock {E}lliptic {G}aussian random processes.
\newblock {\em Rev. Mathem. Iberoamericana}, 13(1):19--89, 1997.

\bibitem{HB}
H.~Bierm\'e.
\newblock {\em Champs aléatoires: autosimilarit\'e, anisotropie et \'etude
  directionnelle}.
\newblock PhD thesis, Universit\'e d'Orl\'eans,
  \url{www.math-info.univ-paris5.fr/~bierme}, 2005.

\bibitem{ABAE}
A.~Bonami and A.~Estrade.
\newblock Anisotropic analysis of some {G}aussian models.
\newblock {\em J. Fourier Anal. Appl.}, 9:215--236, 2003.

\bibitem{Chan}
G.~Chan.
\newblock An effective method for simulating {G}aussian random fields.
\newblock In {\em Proceedings of the statistical Computing section}, pages
  133--138, \url{www.stat.uiowa.edu/~grchan/}, 1999. Amerir. Statist.

\bibitem{Coeurjolly}
J.~F. Coeurjolly.
\newblock {\em Inférence statistique pour les mouvements browniens
  fractionnaires et multifractionnaires}.
\newblock PhD thesis, Université Joseph Fourier, 2000.

\bibitem{Coeurjolly2}
J.~F. Coeurjolly.
\newblock Estimating the parameters of fractional {B}rownian motion by discrete
  variations of its sample paths.
\newblock {\em Stat. Inference Stoch. Process.}, 4:199--227, 2001.

\bibitem{dacunha}
D.~Dacunha-Castelle and M.~Duflo.
\newblock {\em Probabilités et statistiques}, volume~2.
\newblock Masson, 1983.

\bibitem{Dietrich}
C.~R. Dietrich and G.~N. Newsam.
\newblock Fast and exact simulation of stationary gaussian processes through
  circulant embedding of the covariance matrix.
\newblock {\em SIAM J. Sci. Comput.}, 18(4):1088--1107, 1997.

\bibitem{Enriquez}
N.~Enriquez.
\newblock A simple construction of the fractional brownian motion.
\newblock {\em Stochastic Process. Appl.}, 109(2):203--223, 2004.

\bibitem{IL}
J.~Istas and G.~Lang.
\newblock Quadratic variations and estimation of the local {H}ölder index of a
  {G}aussian process.
\newblock {\em Ann. Inst. Henri Poincaré, Prob. Stat.}, 33(4):407--436, 1997.

\bibitem{LESI}
R.~Jennane, R.~Harba, E.~Perrin, A.~Bonami, and A.~Estrade.
\newblock Analyse de champs browniens fractionnaires anisotropes.
\newblock {\em 18eme colloque du GRETSI}, pages 99--102, 2001.

\bibitem{KaplanKuo}
L.~M. Kaplan and C.~C.~J. Kuo.
\newblock An {I}mproved {M}ethod for 2-d {S}elf-{S}imilar {I}mage {S}ynthesis.
\newblock {\em IEEE Trans. Image Process.}, 5(5):754--761, 1996.

\bibitem{Kent}
J.~T. Kent and A.~T.~A. Wood.
\newblock Estimating the fractal dimension of a locally self-similar {G}aussian
  process by using increments.
\newblock {\em J. Roy. Statist. Soc. Ser. B}, 59(3):679--699, 1997.

\bibitem{Lang}
G.~Lang and F.~Roueff.
\newblock Semi-parametric estimation of the {H}\"older exponent of a stationary
  {G}aussian process with minimax rates.
\newblock {\em Stat. Inference Stoch. Process.}, 4(3):283--306, 2001.

\bibitem{Leger}
S.~Leger.
\newblock {\em Analyse stochastique de signaux multi-fractaux et estimations de
  paramètres}.
\newblock PhD thesis, Université d'Orléans,
  \url{http://www.univ-orleans.fr/SCIENCES/MAPMO/publications/leger/these.php},
  2000.

\bibitem{MVN}
B.~B. Mandelbrot and J.~Van~Ness.
\newblock Fractional {B}rownian motion, fractionnal noises and applications.
\newblock {\em Siam Review}, 10:422--437, 1968.

\bibitem{MST}
Y.~Meyer, F.~Sellan, and M.S. Taqqu.
\newblock Wavelets, {G}eneralised {W}hite {N}oise and {F}ractional
  {I}ntegration: {T}he {S}ynthesis of {F}ractional {B}rownian {M}otion.
\newblock {\em J. Fourier Anal. Appl.}, 5(5):465--494, 1999.

\bibitem{Norros}
I.~Norros and P.~Mannersalo.
\newblock Simulation of {F}ractional {B}rownian {M}otion with {C}onditionalized
  {R}andom {M}idpoint {D}isplacement.
\newblock Technical report, Advances in Performance analysis,
  \url{http://vtt.fi/tte/tte21:traffic/rmdmn.ps}, 1999.

\bibitem{Peltier}
R.~F. Peltier and J.~Lévy Véhel.
\newblock Multifractional {B}rownian motion: definition and preliminary
  results.
\newblock Technical report, INRIA, \url{http://www.inria.fr/rrrt/rr-2645.html},
  1996.

\bibitem{LESILeana}
E.~Perrin, R.~Harba, C.~Berzin-Joseph, I.~Iribarren, and A.~Bonami.
\newblock nth-order fractional {B}rownian motion and fractional {G}aussian
  noises.
\newblock {\em IEEE Trans. Sign. Proc.}, 45:1049--1059, 2001.

\bibitem{Ileana}
E.~Perrin, R.~Harba, R.~Jennane, and I.~Iribarren.
\newblock Fast and {E}xact {S}ynthesis for 1-{D} {F}ractional {B}rownian
  {M}otion and {F}ractional {G}aussian {N}oises.
\newblock {\em IEEE Signal Processing Letters}, 9(11):382--384, 2002.

\bibitem{Pipiras}
V.~Pipiras.
\newblock Wavelet-based simulation of fractional {B}rownian motion revisited.
\newblock Preprint, \url{http://www.stat.unc.edu/faculty/pipiras}, 2004.

\bibitem{ramm}
A.~G. Ramm and A.~I. Katsevich.
\newblock {\em The {R}adon {T}ransform and {L}ocal {T}omography}.
\newblock CRC Press, 1996.

\bibitem{Stein}
M.~L. Stein.
\newblock Fast and exact simulation of fractional {B}rownian surfaces.
\newblock {\em J. Comput. Graph. Statist.}, 11(3):587--599, 2002.

\end{thebibliography}
\end{document}